\documentclass{gIPE2e}

\usepackage{epstopdf}
\usepackage{subfigure}
\usepackage[utf8]{inputenc}
\usepackage{array}
\usepackage{subfigure}

\def\R{\mathbb{R}}
\def\L{\mathbb{L}}
\def\H{\mathbb{H}}

\def\<{\langle}
\def\>{\rangle}

\theoremstyle{plain}

\theoremstyle{remark}

\theoremstyle{definition}

\begin{document}


\title{\textit{Optimal parameters identification and sensitivity study for Abrasive Waterjet Milling model}}

\author{
\name{D. Auroux\textsuperscript{a} and V. Groza\textsuperscript{a}$^{\ast}$\thanks{$^\ast$Corresponding author. Email: groza@unice.fr}}
\affil{\textsuperscript{a}Laboratoire de Math{\'e}matiques J.A. Dieudonn{\'e}, Universit{\'e} de Nice - Sophia Antipolis, Nice, France}
}

\maketitle

\begin{abstract}
In this paper we present the work related to the parameters identification for Abrasive Waterjet Milling (AWJM) model that appears as an ill-posed inverse problem. The necessity of studying this problem comes from the industrial milling applications where the possibility to predict and model the final surface with high accuracy is one of the primary tasks in the absence of any knowledge of the model parameters that should be used. The adjoint approach based on corresponding Lagrangian gives the opportunity to find out the unknowns of the AWJM model and their optimal values that could be used to reproduce the required trench profile. Due to the complexity of the nonlinear problem and the large number of the model parameters, we use an automatic differentiation (AD) software tool. This approach also gives us the ability to distribute the research on more complex cases and consider different types of model errors and 3D time dependent model with variations of the jet feed speed. This approach gives us a good opportunity to identify the optimal model parameters and predict the surface profile both with self-generated data and measurements obtained from the real production. Considering different types of model errors allows us to receive the results acceptable in manufacturing and to expect the proper identification of unknowns.
\end{abstract}


\begin{classcode}35B30; 35Q93; 35R30; 49Kxx; 65K10 \end{classcode}

\section{Introduction}

One of the crucial problems in the variety of natural, industrial and economical phenomena which could be modeled by (systems of) partial differential equations (PDEs), is the identification of model parameters by suitable comparison between the experimental measurements referred to the real systems and the predictions of the mathematical models.

In order to solve the direct problem, it is necessary to know in advance all the involved parameters such as coefficients or sources characteristics, but usually different model parameters are unknown or inaccessible and have to be identified only from experimental measurements. However, parameter identification is generally an ill-posed inverse problem \cite{Lavr1, Tarantola}, even if the PDEs are linear under some considerations due to the measurement noises and modeling errors. This kind of problems can be overcome by using so-called regularization methods \cite{Aster, Tikh_Arsenin, Tikh_Glas}, studied by many authors \cite{Tautenhahn2, Denisov, Barbara2, Heinz_Hanke}, and parameters identification could be reformulated in a stable case as a minimization problem with a data mismatch and a regularization term.

In this paper we present the mathematical method to identify unknown parameters of the proposed generic Abrasive Waterjet Milling (AWJM) model which was previously described and studied in \cite{Dragos1, Dragos2, Dragos3, Pablo} and was developped according to the industrial needs for waterjet footprints prediction.

Assuming that model parameters and source terms are known, one can find the shape of the trench profile. In this context this is the forward problem that in our case involves a nonlinear PDE. Some studies of linear  AWJM inverse problems were previously presented in \cite{Bilbao}.

Our goal is to determine the model parameters inaccessible from the experiments to predict the surface construction before performing the manufacturing simulations. Following this requirement we pose the minimization problem, that could be presented as searching the minimum of a cost function -- the difference between the observations and the corresponding model solution. This functional can be minimized using a numerical optimization method based on the limited memory Broyden--Fletcher--Goldfarb--Shanno (BFGS) algorithm \cite{Gil89, Liu89, Didier1} underlying the N2QN1 minimization package from the INRIA MODULOPT library \cite{N2QN1}. This is a gradient descent algorithm which is an iterative process requiring the gradient of the cost function. In our work the gradient vector is obtained numerically using the automatic differentiation software TAPENADE \cite{TAPENADE}, which can be interpreted as the Lagrange multiplier of the model equations in terms of adjoint problem \cite{Quarteroni, Lions}.

By adding noise to the artificial data, we show that in fact the parameter identification problem is highly unstable and strictly depends on input measurements. Furthermore, we demonstrate that Tikhonov regularization could be effectively used to deal with the presence of data noise and to improve the identification correctness.

The paper is organized as follows. In the first section we present the proposed mathematical model of the problem, explain the adjoint approach based on corresponding Lagrangian that was used to get the gradient of the cost function, needed in the implementation of the gradient descent L-BFGS based algorithm. We also provide a short list of gradient descent algorithms with some explanations. Section \ref{sec_Lin2D} represents the parameters identification for linearized 2D cases, based on self-generated input data as well as on real experiment measurements. Section \ref{sec_sensitivity} contains the results corresponding to the sensitivity study of the model in case of various measurement errors and section \ref{sec_concl} finalizes this paper with some conclusions and prospectives of future research.

\section{Mathematical model}
\label{sec_model}  
\subsection{Proposed model}

The model of jet footprints is defined as a nonlinear partial differential equation that describes the forming of the surface for any jet feed speed and milling direction regardless of the target material types.
We consider the time interval $[0,T]$ and denote by $\Omega$ a bounded domain of $\R^2$.

The Abrasive Waterjet Milling model introduced in \cite{Dragos1, Dragos2, Dragos3} is described as:

\begin{equation}\label{model}
\frac{\upartial \bm Z}{\upartial t}=\frac{\bm E(x,y){\rm e}^{a \bm Z}}{\left( 1 + \left|\nabla \bm Z\right|^2 \right)^{k/2}} \qquad\qquad \text{in} \: \Omega \times [0,T], 
\end{equation}
where 
\begin{itemize}
\item $(x,y) \in \Omega$
\item$\bm Z(x,y,t) \in \H^1_0(\Omega,\R^+)$ is the parametrization of the surface,
\item$a,k \in \R^+$ are the model parameters, 
\item$\bm E(x,y) \in \L^2(\Omega)$ is the Etching rate function which is also model parameter,
\end{itemize}
with initial and boundary conditions:
$$
\begin{cases}
\bm Z(x,y) = 0 \qquad \text{\quad on} \quad \upartial\Omega,\\
\bm Z(x,y) = \bm Z_0 \qquad\; \text{at} \quad t=0 .
\end{cases}
$$

Further the AWJM model could be rewritten in general as:
\begin{equation}\label{general model}
\frac{\upartial \bm Z}{\upartial t}=\bm F(\bm Z, t,\bm u),
\end{equation}
where $\bm u=\{a,k,\bm E\}$ is a set of model parameters.


The proposed model above explains the impact of an abrasive waterjet of radius $a$ on the flat surface at the initial time moment, which changes its shape during the process and movement of the jet in the $y$-direction under the action of Etching rate function $\bm E(x,y)$. The etching rate function, waterjets radius and change of the trenches profile accordingly to the depth and distance between the surface and the jets nozzle have influence on the intensity of the jet impact, and on the modification of the surface, which is described here as a function $\bm Z(x,y,t)$. To model more specifically some of these processes, the exponential component and the gradient of the surface parametrization were introduced. Due to the dependency of Etching rate function $\bm E$ on the properties of the workpiece material and machine parameters as pressure, velocity and abrasives mass flow, it should be calibrated using the experimental data.

\begin{figure}
     \begin{center}
        \subfigure[Cross section]{%
            \label{jet_cross}
				\resizebox*{6cm}{!}{\includegraphics{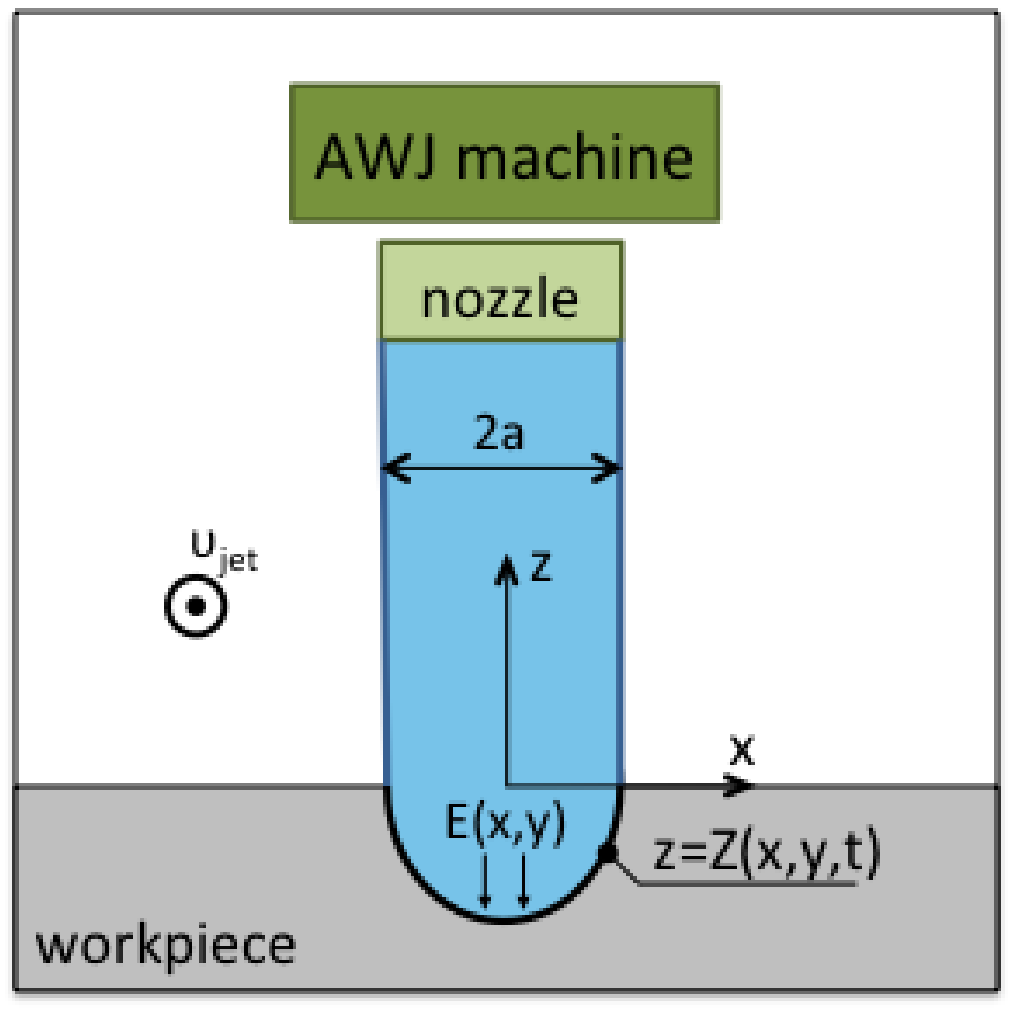}}}\hspace{5pt}
        \subfigure[Top view]{%
           	\label{jet_top}
	           	\resizebox*{6cm}{!}{\includegraphics{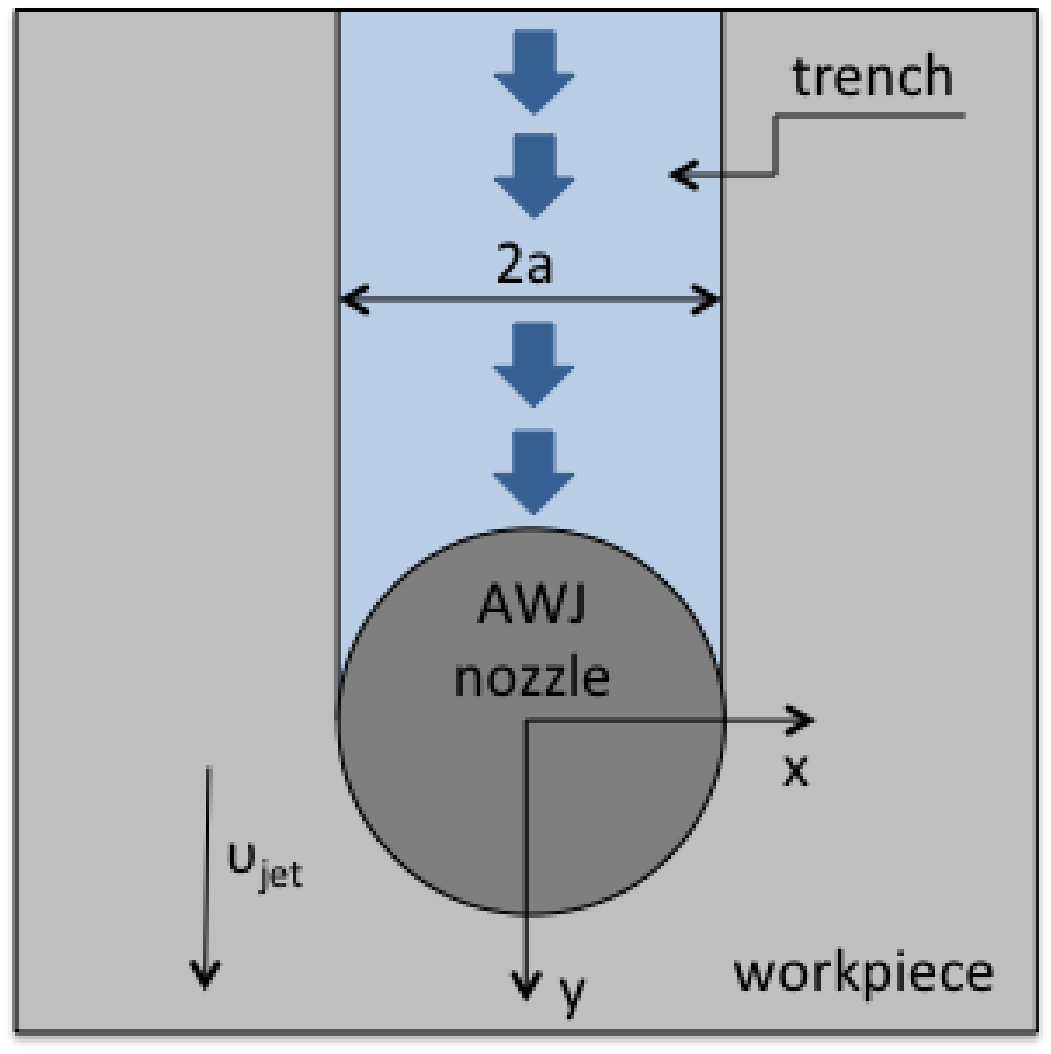}}}
    \caption{Schematic of the AWJM process and jet footprint.
     }%
   \label{Schema}
   \end{center}
\end{figure}

Schematically AWJM process and jet footprint at some time moment $t \in [0,T]$ could be presented as on the Figure \ref{Schema}.

\subsection{Cost function and adjoint model}

For a given set of parameters $\bm u=\{a,k,\bm E\}$ we formulate a minimization problem in order to find $\bm u^*$ such that
\begin{equation}\label{minprob}
\bm J(\bm u^*)= \inf_{\bm u}\bm J(\bm u)  ,
\end{equation}
where $\bm J(\bm u)$ is the cost function measuring the difference between the model solution and experimental observations
\begin{equation}\label{general cost}
\bm J(\bm u)=\int\limits_\Omega \|\bm Z(x,y,T)-\bm Z_{\rm exp}(x,y)\|^2 {\rm d}x{\rm d}y + \alpha\|\bm u-\bm u_{\rm b}\|^2 ,
\end{equation}
under the constraint that $\bm Z$ is the solution of our model with inputs $\bm u \in \R^+ \times \R^+ \times \L^2(\Omega,\R^+)$.

Here $\bm Z_{\rm exp}$ are the experimental measurements, $\bm u_{\rm b}$ is an a priori estimation of the set of parameters $\bm u$, and $\alpha>0$ is the Tikhonov regularization coefficient, that plays an important role in bringing the adjoint problem to a well-posed form and helping to identify suitable unknown parameters.

We have already introduced the solution of the parameters identification problem as a minimization problem of the mismatch between solution and experiment measurements, and in order to overcome the instability and reduce the inaccuracy, regularization theory suggests to augment the data mismatch by a regularization term, which is known as Tikhonov regularization \cite{Tikh_Arsenin}.

In order to find an optimal $\bm u^*$ one needs the gradient of the cost function $\nabla \bm J(\bm u)$ for the minimization process based on iterative gradient descent principles as quasi-Newton algorithms.

The minimization problem \eqref{minprob} could be regarded as an optimal control problem, and the approach based on Lagrangian multipliers could be used to get the solution of it \cite{Quarteroni, Lions}. The solution of the optimal control problem can therefore be regarded as the search of critical "points" of the Lagrangian functional $\bm L(\bm u,\bm Z,\bm P)$ associated to the constrained minimization, that we introduce as following:
\begin{equation}\label{Lagrangian}
\bm L(\bm u,\bm Z,\bm P) = \bm J(\bm u) + \int\limits_0^T \int\limits_\Omega \bm P(t,x,y)\left(\frac{\upartial \bm Z}{\upartial t} - \bm F(\bm Z, t,\bm u)\right){\rm d}x{\rm d}y{\rm d}t,
\end{equation}
where $\bm P$ is the Lagrangian multiplier associated to the constraint that $\bm Z$ is a solution of \eqref{model} from the same space $\H^1_0(\Omega)$. 

We get the optimal system
\begin{eqnarray}
&& \label{partP} \frac{\upartial \bm L}{\upartial \bm P} = 0, \\
&& \label{partZ} \frac{\upartial \bm L}{\upartial \bm Z} = 0, \\
&& \label{partu} \frac{\upartial \bm L}{\upartial \bm u} = 0,
\end{eqnarray}
from where the statement \eqref{partP} says that $\frac{\upartial \bm Z}{\upartial t} = \bm F(\bm Z, t,\bm u)$, where $\bm Z$ is solution of the direct model. We can obtain the adjoint equation, where $\bm P$ should be also a solution from the expression \eqref{partZ}. And from \eqref{partu} we have the gradient of the cost function from the initial adjoint state $\nabla \bm J(\bm u) = \bm P(0)$, if $\bm Z$ is a solution of the direct problem \eqref{general model}, corresponding to the set of parameters $\bm u$.

From equations \eqref{Lagrangian} and \eqref{partZ} we get by integration by parts:
\begin{equation}
\frac{\upartial \bm P}{\upartial t} =  -\left( \frac{\upartial \bm F}{\upartial \bm Z} \right)^{\rm T} \bm P + \left( \bm Z - \bm Z_{\rm exp}\right)\bm \delta(t-T), \\
\end{equation}

\begin{align*}
\bm P(x,y,T) = 0, \\
\bm P(x,y,t)\mid_{\upartial\Omega}\: = 0.
\end{align*}

To find out the optimal model parameters we have to consider the discrete system, not the continuous one. Indeed we need to minimize the discrete cost function and it requires the gradient of the discrete cost function. For that the discretized adjoint statement could be obtained from the disretized forms of the Lagrangian and its discrete derivatives, but more efficient way is to use automatic differentiation software (i.e. TAPENADE) which bases on the principles explained above. It computes and provides the gradient of the discretized cost function which is further used to solve the minimization problem \eqref{minprob}. Then all the optimal model parameters can be computed from only the trench profile, in other words from one experiment done on AWJ machine.

\subsection{Minimization process}

The most common minimization techniques are gradient descent algorithms which are iterative processes converging to the minimum. The minimization of the cost function could be realized in such a way, where the descent step is computed by quasi-Newton type algorithms. Particularly we used the L-BFGS based algorithms implemented in N2QN1 minimization package from "MODULOPT" library \cite{N2QN1}.

Starting from the initial given position $\bm u_0$ (which is defined from some background estimation $\bm u_0 = \bm u_{\rm b}$), one has to compute the cost function $\bm J(\bm u_j)$ on each iteration and its gradient $\nabla \bm J(\bm u_j)$ to shift to the next discretized step and update the solution:
$$
\bm u_{j+1} = \bm u_j - \bm \rho_j \nabla \bm J(\bm u_j) ,
$$
where $\bm \rho_j > 0$ is a descent step of minimization, and $j > 0$ is the iteration number.

This is a first order scheme which could be justified from the first order Taylor expansion of the cost function $\bm J(\bm u)$ for any two iterations $\bm u_j$ and $\bm u_{j+1}$:
$$
\bm J(\bm u_{j+1}) = \bm J(\bm u_j - \bm \rho_j \nabla \bm J(\bm u_j)) \simeq \bm J(\bm u_j) - \bm \rho_j \| \nabla \bm J(\bm u_j)\|^2 \le \bm J(\bm u_j).
$$

The classical Newton type of minimization algorithms uses $- \bm \rho_j \bm H_j^{-1} \nabla \bm J(\bm u_j)$ as the direction of descent, where $\bm H_j=\nabla^2 \bm J(\bm u_j)$ is the Hessian of the cost function $\bm J$ at any iteration $j$. Due to the size of the problem one can not compute the Hessian directly nor its inverse.

In quasi-Newton algorithms the Hessian does not need to be computed and inverted, but it can be replaced by a symmetric positive approximation $(\bm Q_j)$ to $\bm H^{-1}$.

The L-BFGS \cite{Liu89} is a BFGS algorithm in the quasi-Newton's family methods which uses an approximation of the inverse Hessian, but with limited amount of memory: it stores only $m$ last values of $\bm u_j$ and $\nabla \bm J(\bm u_j)$. The update formulas for L-BFGS algoritm are presented in \cite{Gil89}.

\section{Model parameters identification}
\label{sec_Lin2D}

\subsection{Identification based on self-generated 2D profiles}
\label{subsec_Lin2D}

From equation \eqref{model} we get the AWJM model corresponding to steady problems for non-moving jets, with initial and boundary conditions:
\begin{equation}\label{2Dmodel}
\frac{\upartial \bm Z}{\upartial t}=\frac{\bm E(x){\rm e}^{a\bm Z}}{\left( 1 + {\bm Z_{x}}^2 \right)^{k/2}} \;,
\end{equation}
\begin{align*}
\bm Z\mid_{t=0}\: = 0,\\
\bm Z\mid_{\upartial\Omega_1}\: = 0.
\end{align*}

In this case we define a symmetric domain $\Omega_1 = \{ x : x\in [-x_1;x_1 ] \} $, where $x_1$ always depends on the actual experimental parameters or measurements, due to the changes of the jets radius.
Parameter identification problems usually bring many difficulties due to model errors and measurement noise. First we assume that measurement errors $\varepsilon_{\rm exp}$, which come from the noisy measurements, are random variables with a Gaussian probability density function and a zero mean.
To generate the noisy trench profile we introduce the following model:
\begin{equation}\label{2Dmodel_noise}
\frac{\upartial \bm Z}{\upartial t}=\frac{\bm E(x){\rm e}^{a\bm Z}}{\left( 1 + {\bm Z_{x}}^2 \right)^{k/2}} + \lambda\,\bm \varepsilon_{\rm exp} \;,
\end{equation}
where $\lambda$ is the factor corresponding to the percentage of the applied calibrated uncorrelated noise, $\bm \varepsilon_{\rm exp} \sim \mathcal N(0,1)$.

\subsection{Discretization parameters for self-generated data}
\label{subsec_Discr2D}

For the first numerical experiments, "pseudo-experimental" surface profiles are generated with arbitrary values of model parameters $a=2$, $k=3$ and Etching rate function $\bm E_0(x)$ defined on $\Omega_1 = \{ x : x\in [-1;1 ] \} $. To discretize the domain $\Omega_1$, a regular grid of 200 points with a step $\Delta x = 0.01$ is used and the time period is taken as unit $t\in\left[ 0,1 \right]$ with $\Delta t = \frac{\Delta x^2}{4}$. To solve our minimization problem numerically the central difference method is used for the time integration.

The identification of both model parameters $a$ and $k$ is realized in the conditions of noise applied with different levels (10 and 30\%) to the input surface. A Tikhonov regularization term is used with different values of the factor $\alpha$. This allows us to neutralize problems in identification process regarding the ill-posedness and to accelerate the minimization process. Using various values of $\alpha$, we show the importance of the regularization. 

\begin{figure}
     \begin{center}
     \resizebox*{6cm}{!}{\includegraphics{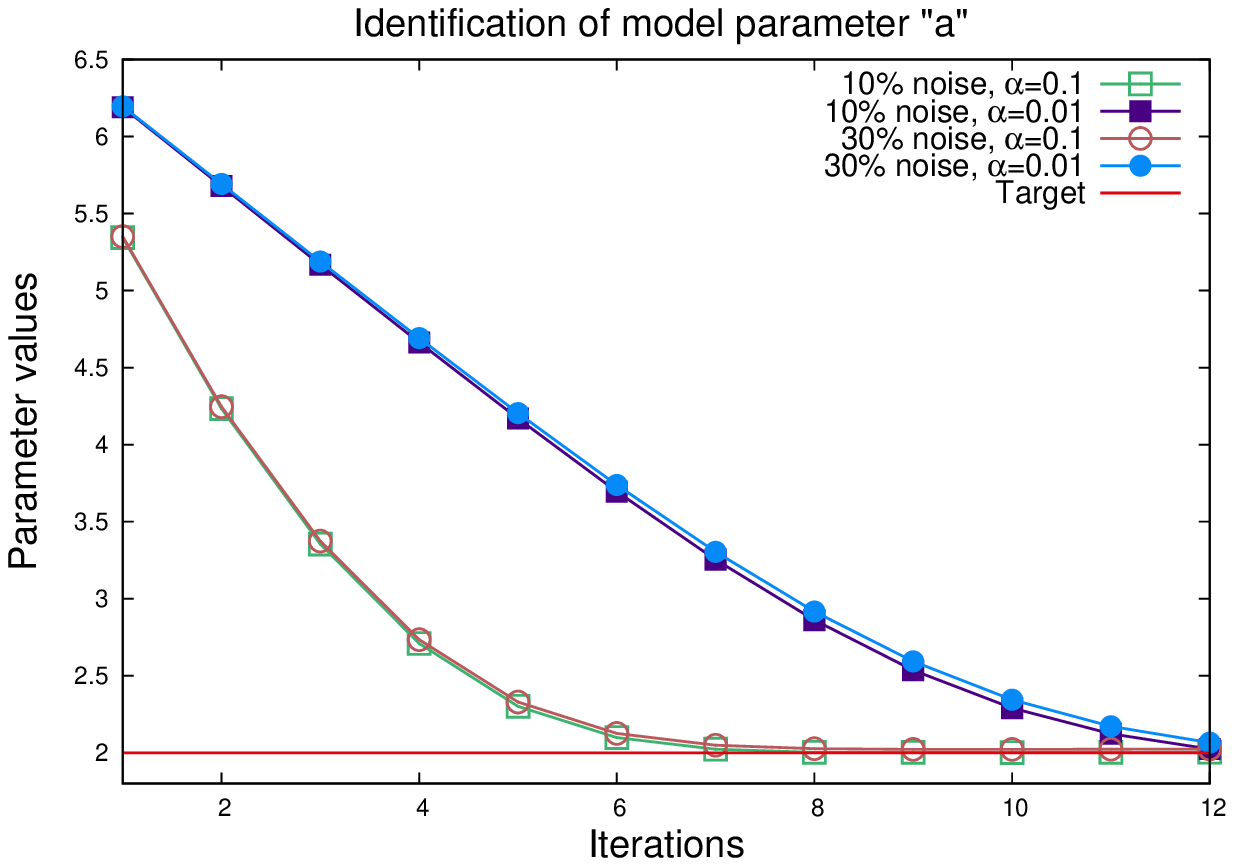}}\hspace{5pt}
     \resizebox*{6cm}{!}{\includegraphics{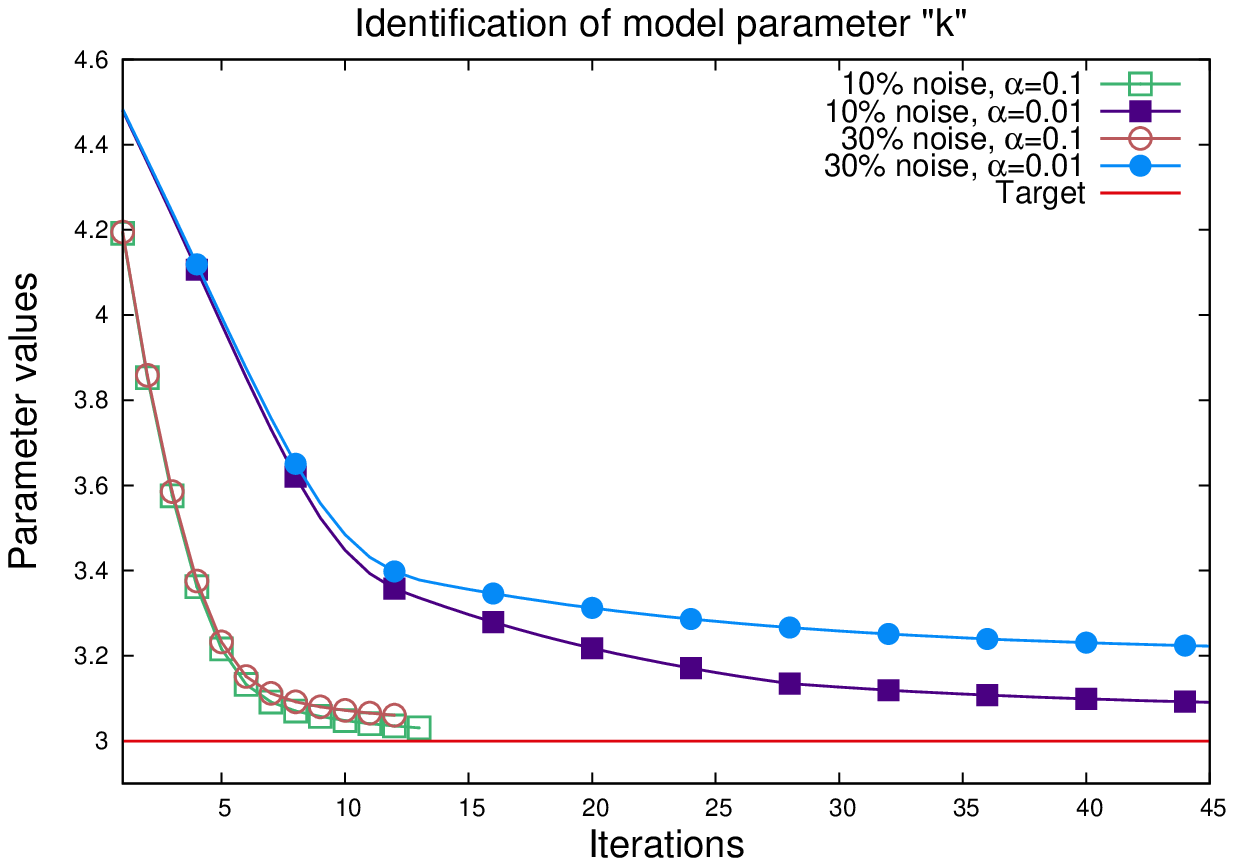}}
    \caption{Results of the identification of model parameters $a$ and $k$ under the conditions of applied noise with levels of $10\%$ and $30\%$, using different Tikhonov regularization coefficients.
     }%
   \label{graphsak}
   \end{center}
\end{figure}

Results and comparisons are shown on Figure \ref{graphsak}. We can notice that the background estimation of the unknown model parameters accelerates the minimization process, but the higher the level of noise in the input data the bigger errors in the identification. Even with a high level of noise, identification works properly and allows us to determine the unknowns being far enough from the first estimations. The choice of the regularization coefficient plays a very important role in the identification problem.

For the identificaiton of AWJM model parameters $a$ and $k$, one can see on Figure \ref{graphsak} that the use of right regularization factors strongly improves the speed and accuracy of the identification also with high levels of noise. In this particular case $\alpha=0.1$ gives the opportunity to identify initial model parameters within 8 and 15 iterations respectively, but with the value $\alpha=0.01$ it takes more time to approach the required precision. The use of too large or small values may prevent to identify the parameters at all.


The first part of the experiments only includes the identification process corresponding to the standard cost function \eqref{general cost}, but in order to obtain a smooth solution ($E(x,y) \in \H^1(\Omega)$ instead of $\L^2(\Omega)$), we now change the regularization term to the one with the gradient of the Etching rate function:

\begin{equation}\label{gradient cost}
\bm J(\bm u)=\int\limits_\Omega \|\bm Z(x,T)-\bm Z_{\rm exp}(x)\|^2 {\rm d}x + \alpha\|\nabla \bm E\|^2 .
\end{equation}

We also need to change the configuration of the minimization process as well, and note that in order to identify the function $\bm E$, another value of Tikhonov regularization term $\alpha = 10^{-6}$ is now the closest to the optimal one. It was obtained by L-curve method which was first applied by Lawson and Hanson \cite{Laws_Hanson} and more recently by Hensen and O'Leary \cite{Hensen_Oleary}.

\begin{figure}
     \begin{center}
        \subfigure[First estimation]{%
            \label{graphsE_a}
				\resizebox*{6cm}{!}{\includegraphics{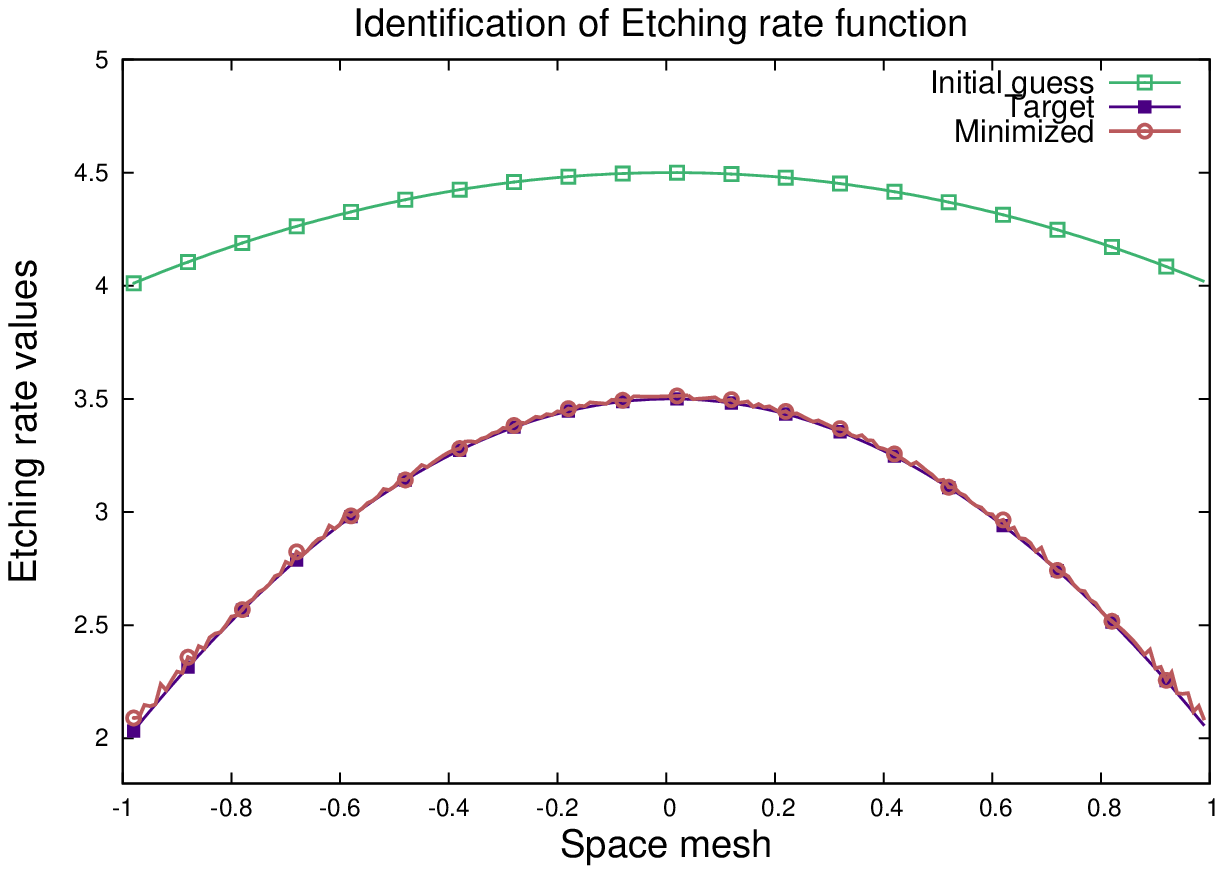}}}\hspace{5pt}
        \subfigure[Second estimation]{%
           	\label{graphsE_b}
	           	\resizebox*{6cm}{!}{\includegraphics{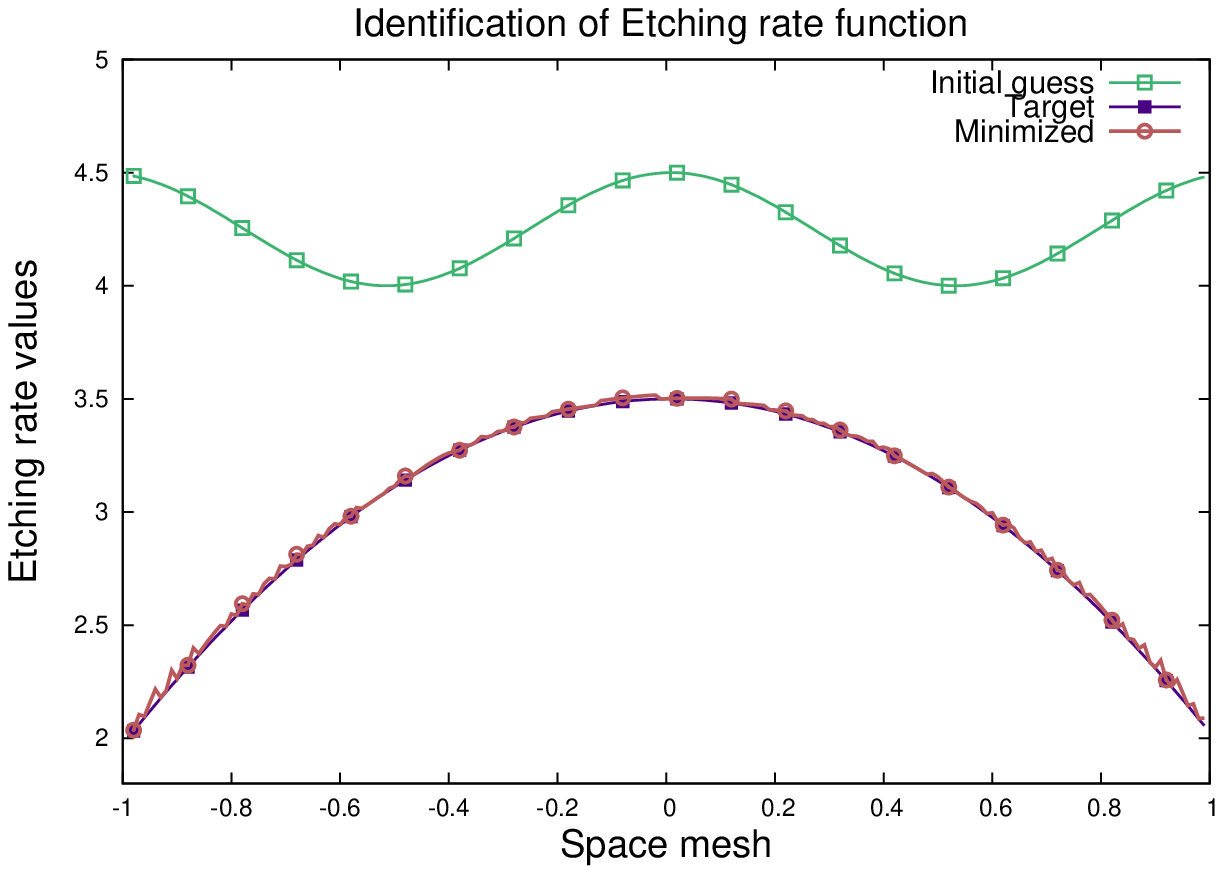}}}
    \caption{Results of the identification of Etching rate function $E$ starting from different initial conditions.
     }%
   \label{graphsE}
   \end{center}
\end{figure}

The goal of the next simulations is to correctly identify the Etching rate function $\bm E$. It was used to generate the surface profile in the perfect conditions without any noises or errors accounted in the model. The other values of the model parameters $a=2, k=3$ were fixed. Results of the identification are presented on Figure \ref{graphsE}. The point of identifying the Etching rate function is to focus on the most influential term of AWJM model, that increases the number of unknowns by several orders.

The results shown on Figure \ref{graphsE} confirm the possibility to identify the Etching rate function in the ideal situation even with a wrong first estimation. Here the total number of unknowns is 200 which is the number of mesh points. In case of a really bad first estimation (Figure \ref{graphsE_b}) the minimization process takes more or less the same time and stops after 19 iterations by reaching the defined accuracy of $10^{-4}$, when in the first case, shown on Figure \ref{graphsE_a}, it takes 16 iterations.

More interesting is the case when the Etching rate function is not so smooth and simple, or evenmore is unknown, and the trench profile has some specialties caused ratherish by measurement errors.

\subsection{Real experiment measurements}
\label{subsec_2Dreal}

Unlike the previous cases, the trench profile is now the only input for our problem. We present the results of identification of the Etching rate function which should be used in the direct AWJM model to be able to reproduce the required trench profile. These measurements were obtained from experiments done in collaboration with STEEP project partners. This data corresponds to experiments with a jet feed speed of 2000 mm/min, and noozle diameter of AWJ machine of $0.5$. We use this value as a background estimation of model parameter $a$. The number of mesh points in this case equals 228 with the step $\Delta x = 0.0048$ and $\Omega_1 = [-0.55;0.55 ] $. The minimization process is realized using N2QN1 minimizer from the "MODULOPT library" with the regularization factor $\alpha=10^{-5}$ that was obtained numerically by L-curve method. In lack of any knowledge about the Etching rate function we started the identification from the zero assumption $\bm E_0 = \bm 0$.

\begin{figure}
     \begin{center}
        \subfigure[Cost function]{%
            \label{cost_behav}
				\resizebox*{6cm}{!}{\includegraphics{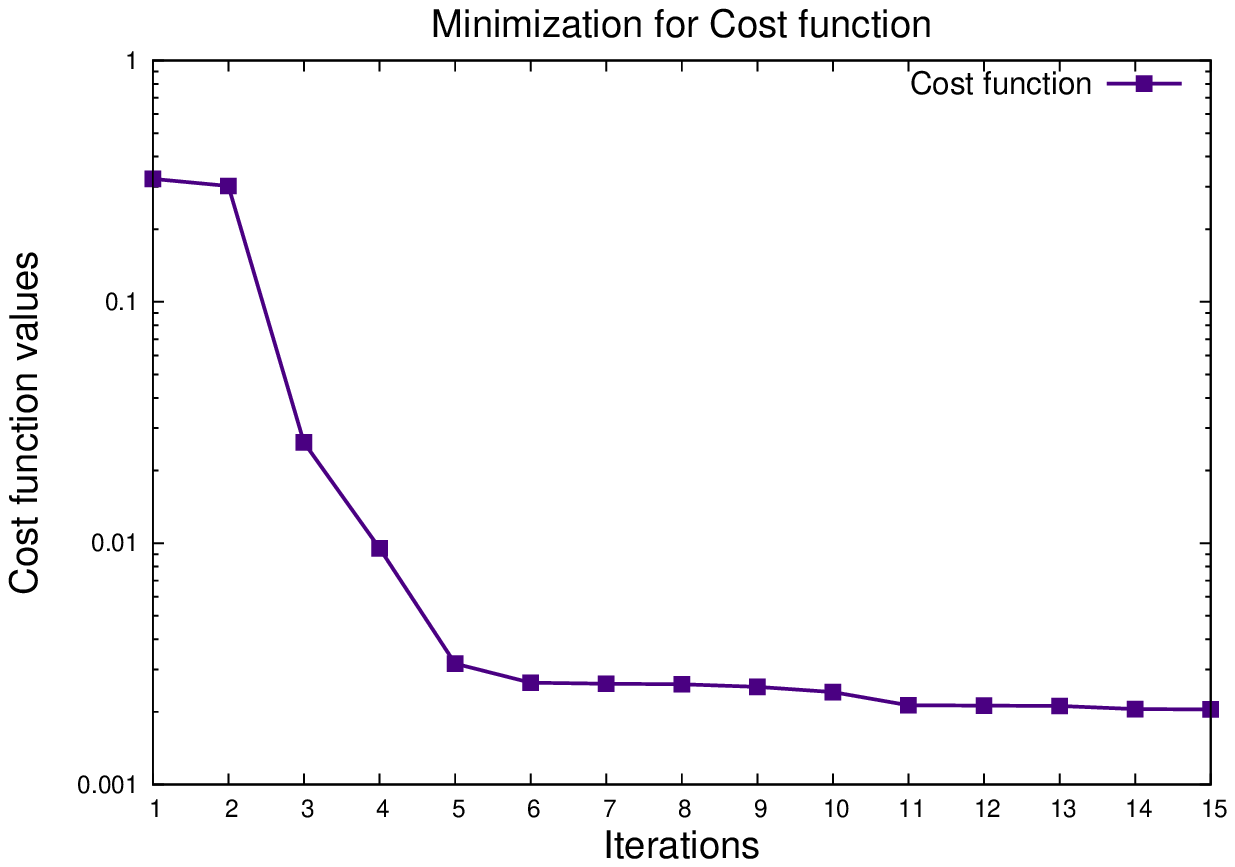}}}\hspace{5pt}
        \subfigure[Gradient of cost function]{%
           	\label{grad_behav}
	           	\resizebox*{6cm}{!}{\includegraphics{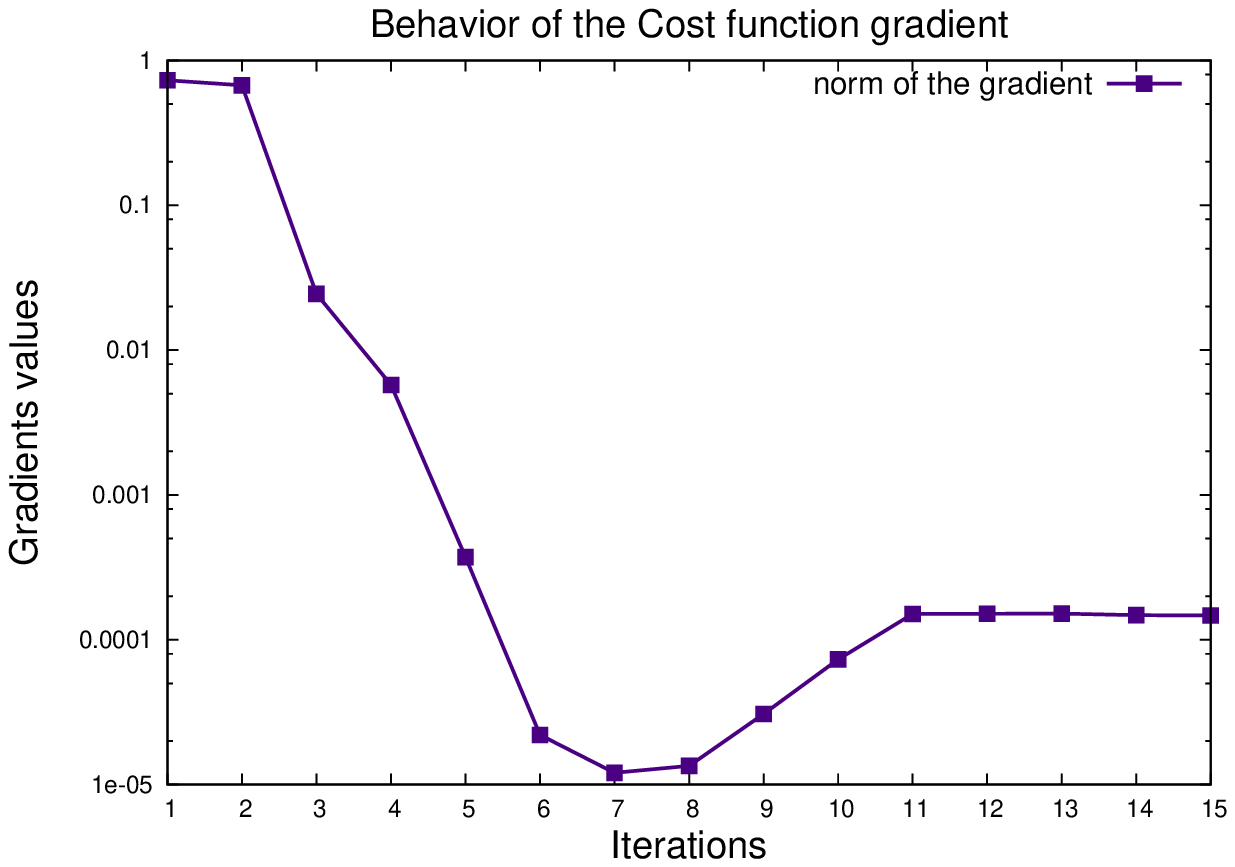}}}
    \caption{Behavior of the cost function and its gradient in the minimization process, corresponding to the real measurement data.
     }%
   \label{behav}
   \end{center}
\end{figure}

Figures \ref{cost_behav} and \ref{grad_behav} show respectively the behavior of the cost function and of the norm of its gradient during the minimization process. It shows that the minimizer N2QN1 works well. Using 15 minimization iterations, the Etching rate function is identified with the accuracy of $4 \times 10^{-3}$.

\begin{figure}
     \begin{center}
     \resizebox*{6cm}{!}{\includegraphics{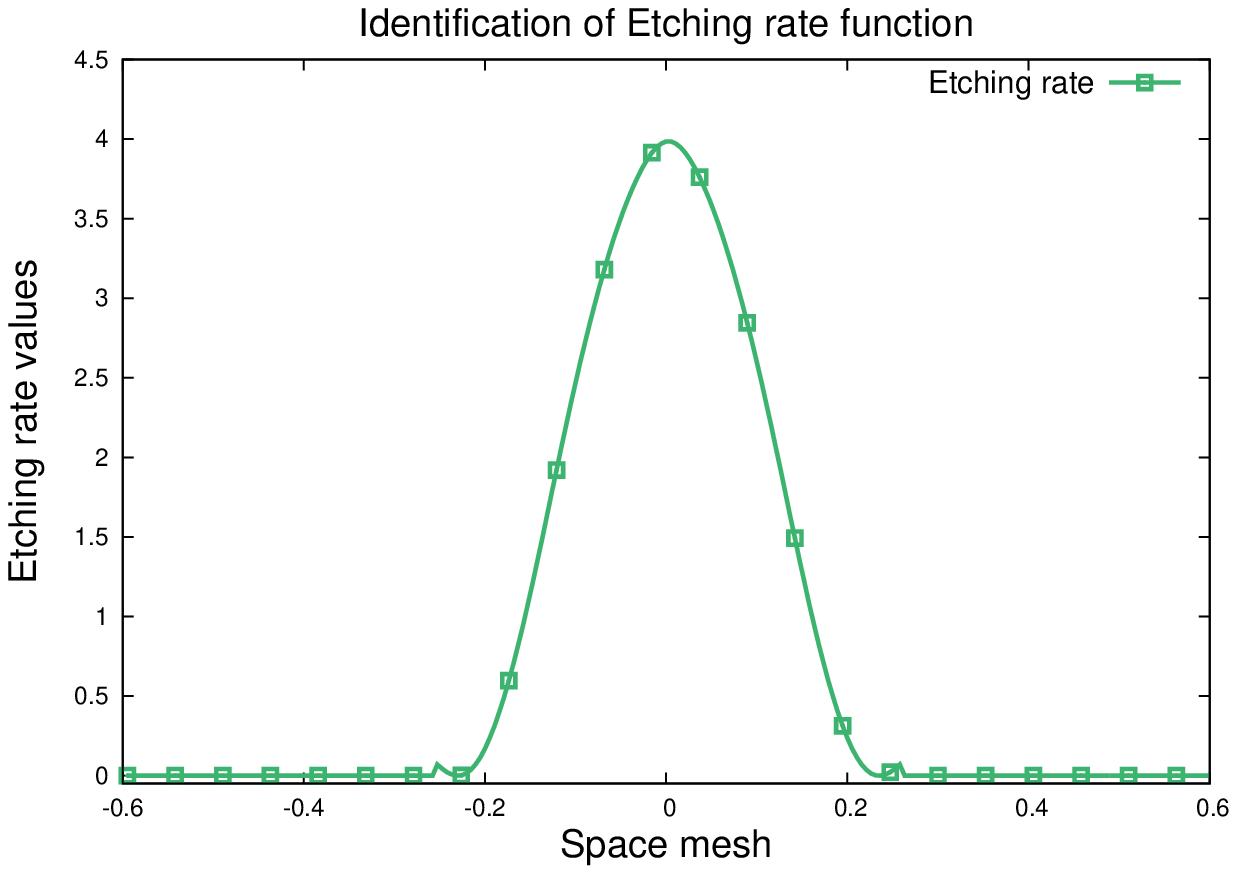}}\hspace{5pt}
     \resizebox*{6cm}{!}{\includegraphics{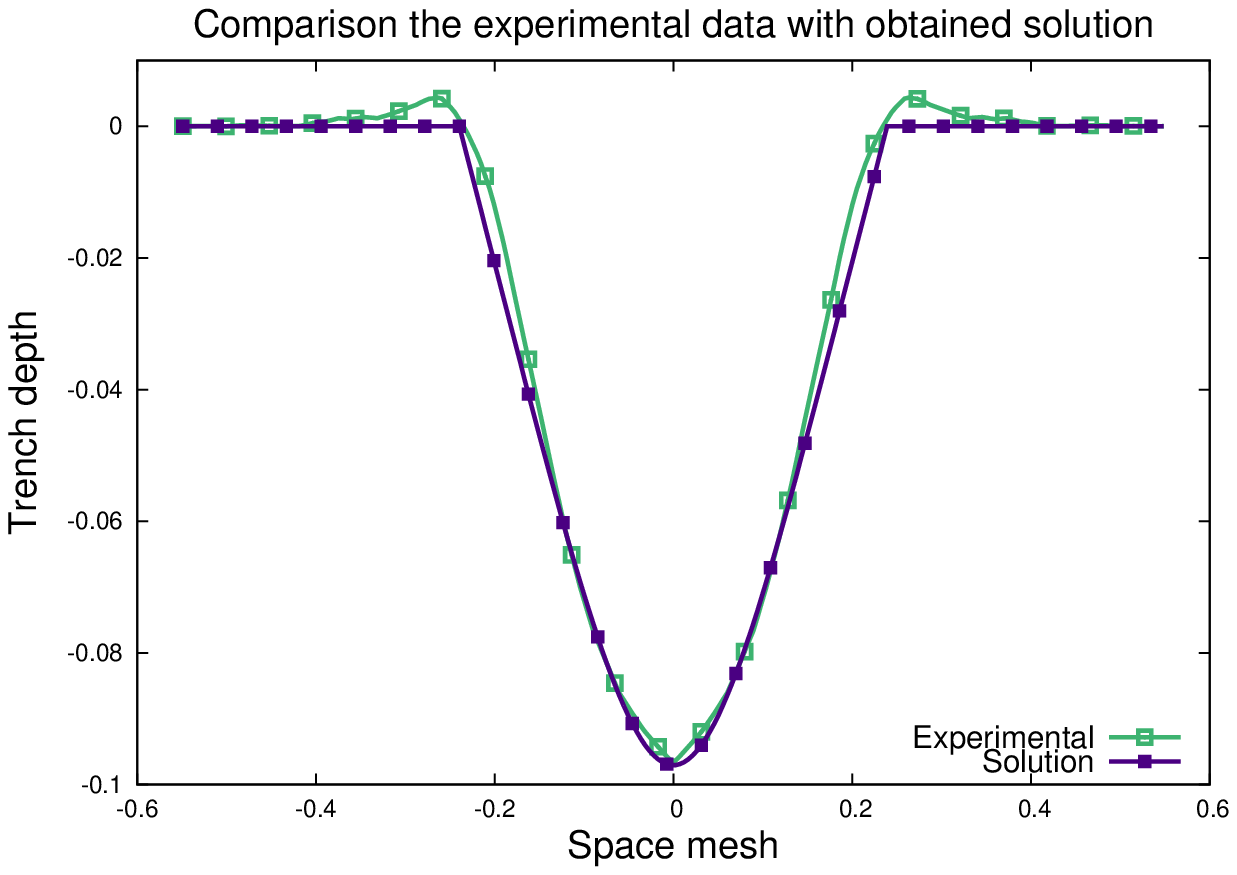}}
    \caption{Results of the identification of Etching rate function and corresponding trench profile in comparison with experimental measurements.
     }%
   \label{realE}
   \end{center}
\end{figure}

The identified Etching rate function from the experimental data allows us to reproduce numerically the trench profile with an error in terms of $L^2$ norm smaller than 3\% (Figure \ref{realE}). The used mathematical model \eqref{model} has some limitations itself and can not cover some effects on the edges of the trench which in practice is explained by redeposition of the surface material. These secondary effects appearing as the result of high power of the waterjet impact should be studied separately and are not considered in our work.

\section{Numerical results of sensitivity study}
\label{sec_sensitivity}

One of the most important things in the identification problems is the sensitivity study which could give us the opportunity to provide correct AWJM model parameters to reproduce required shape of the trench. The point of this section is to study the capacity of the approach, and to observe the possibilities to identify model parameters even with high level of measurement errors which are always present in the provided input data. Based on AWJM model \eqref{2Dmodel_noise} we generate different trench profiles with predefined values of $\bm E=\bm E_0, a=a_0, k=k_0$ respectively to different levels of uncorrelated noise up to 40\%. We use these noisy profiles as the only input of the method and we try to identify the Etching rate function that could satisfy the requirements. As this part of the study is based again on self-generated input data, we use the same numerical parameters and initial assumptions as in Section \ref{subsec_2Dreal}. We also assume that initial Etching rate function $\bm E_0$ has symmetrical gaps corresponding to the edges of the trench. 

Tikhonov regularization has a very strong influence. It allows us to vary the identification results between accuracy and smoothness. We can then choose the optimal values that can be acceptable and suitable in the real experiments and lead to reconstruction of the surface profile with high precision at the same time.

\begin{figure}
     \begin{center}
        \subfigure[Single profile, 5\% of noise]{%
            \label{1tr_5}
				\resizebox*{6cm}{!}{\includegraphics{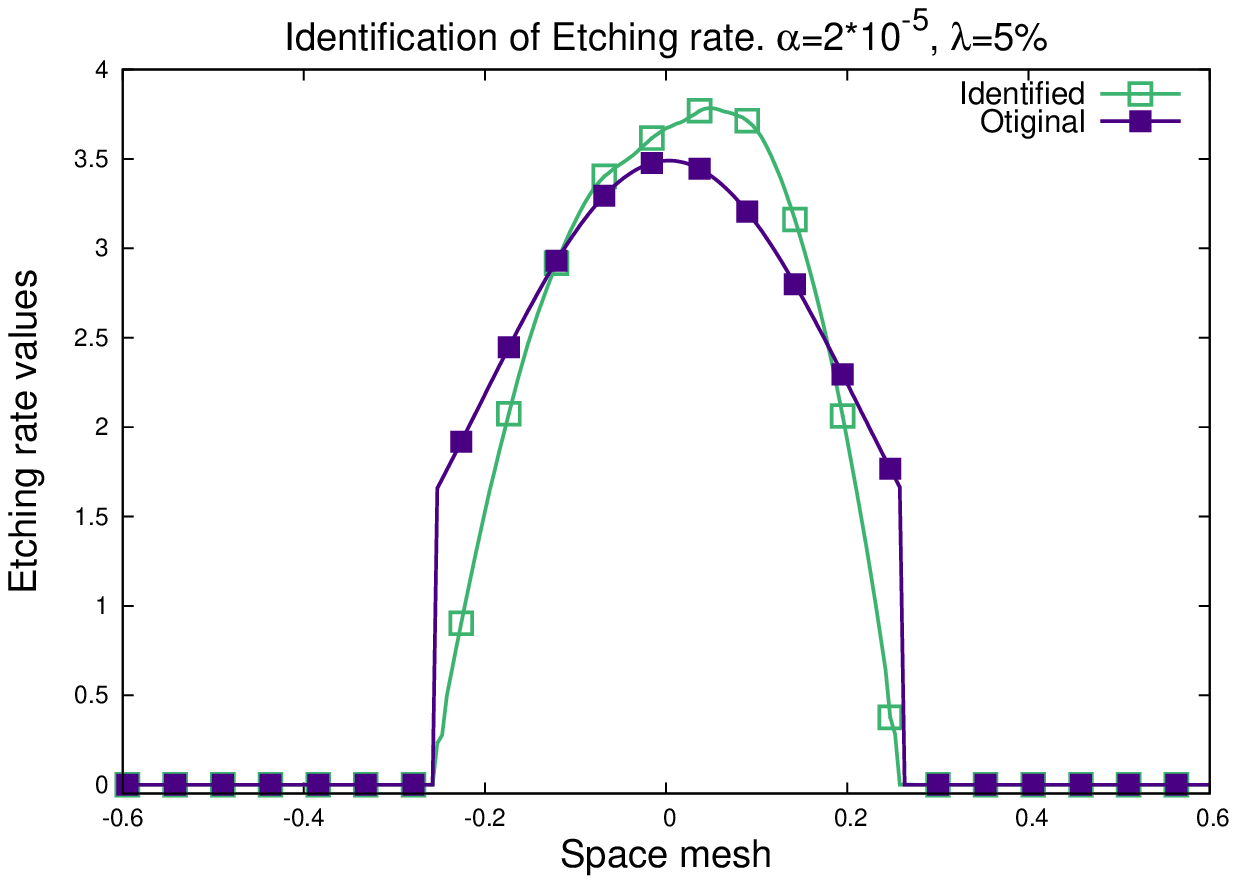}}\hspace{5pt}
				\resizebox*{6cm}{!}{\includegraphics{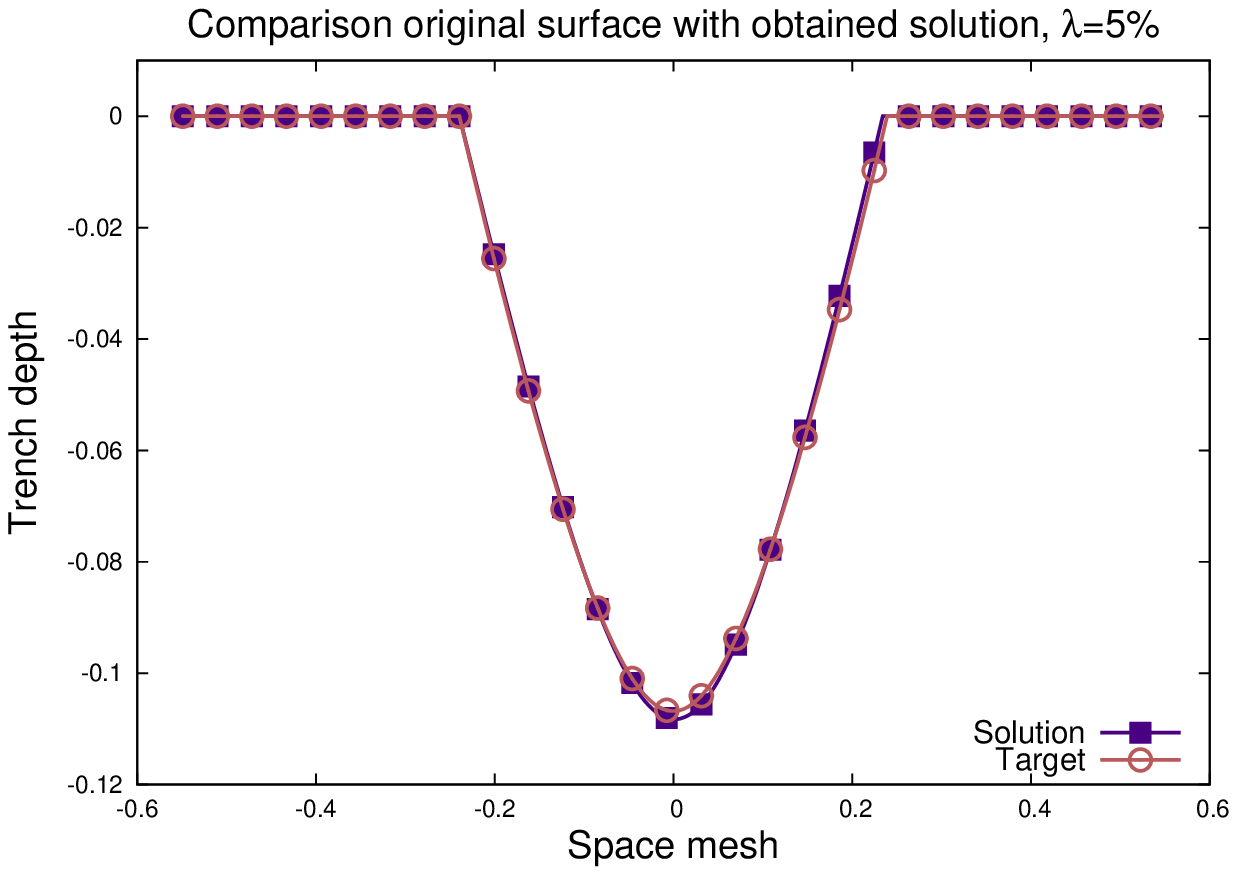}}} \\
        \subfigure[Single profile, 15\% of noise]{%
            \label{1tr_15}
				\resizebox*{6cm}{!}{\includegraphics{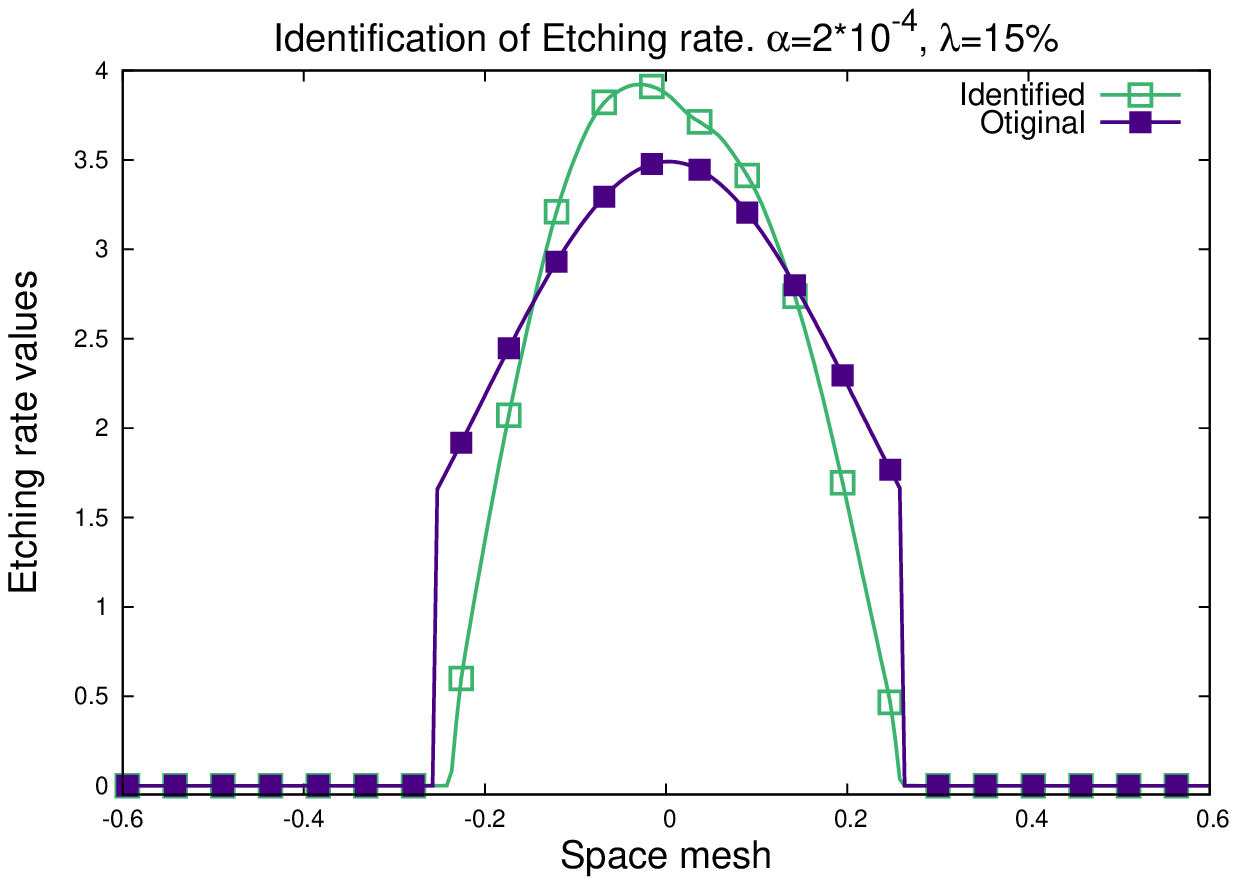}}\hspace{5pt}
				\resizebox*{6cm}{!}{\includegraphics{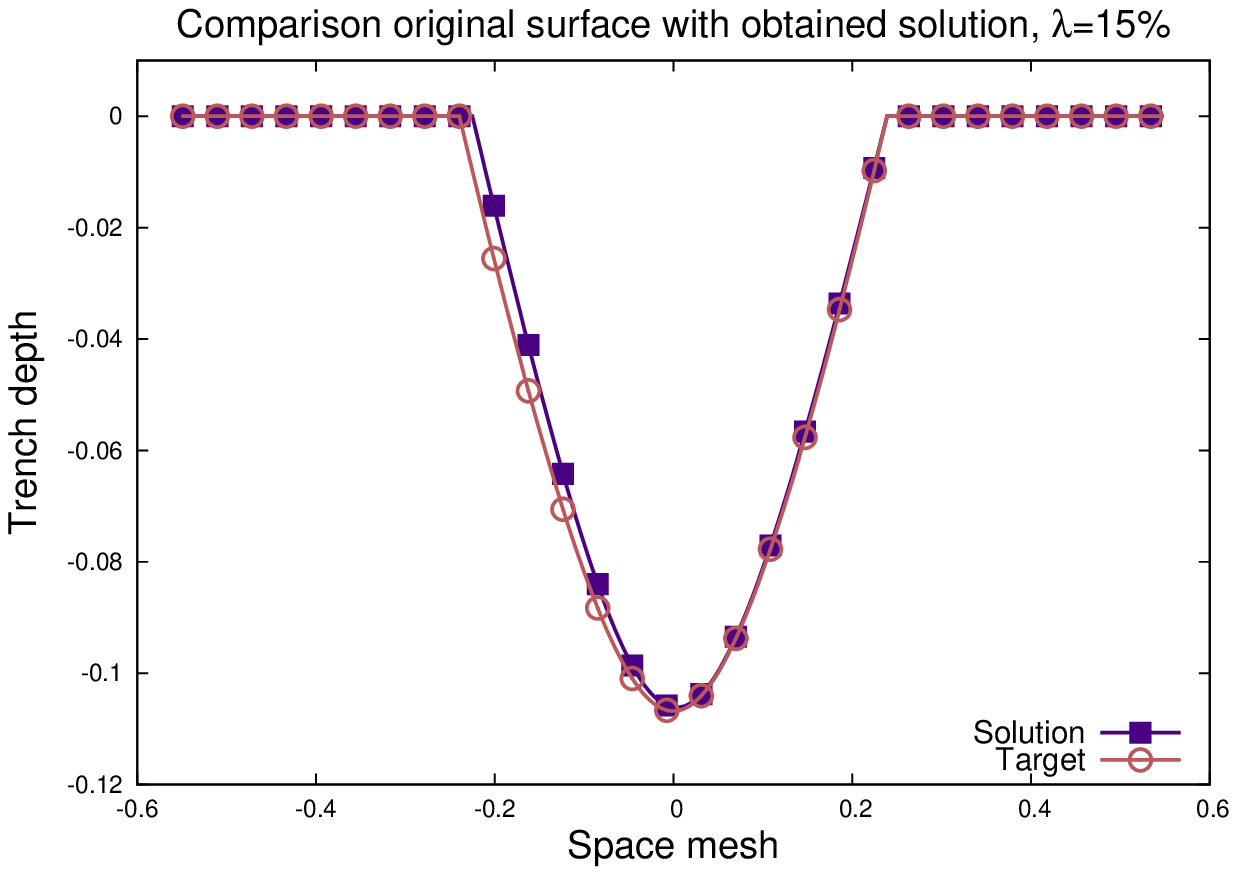}}} \\
		\subfigure[Single profile, 30\% of noise]{%
            \label{1tr_30}
				\resizebox*{6cm}{!}{\includegraphics{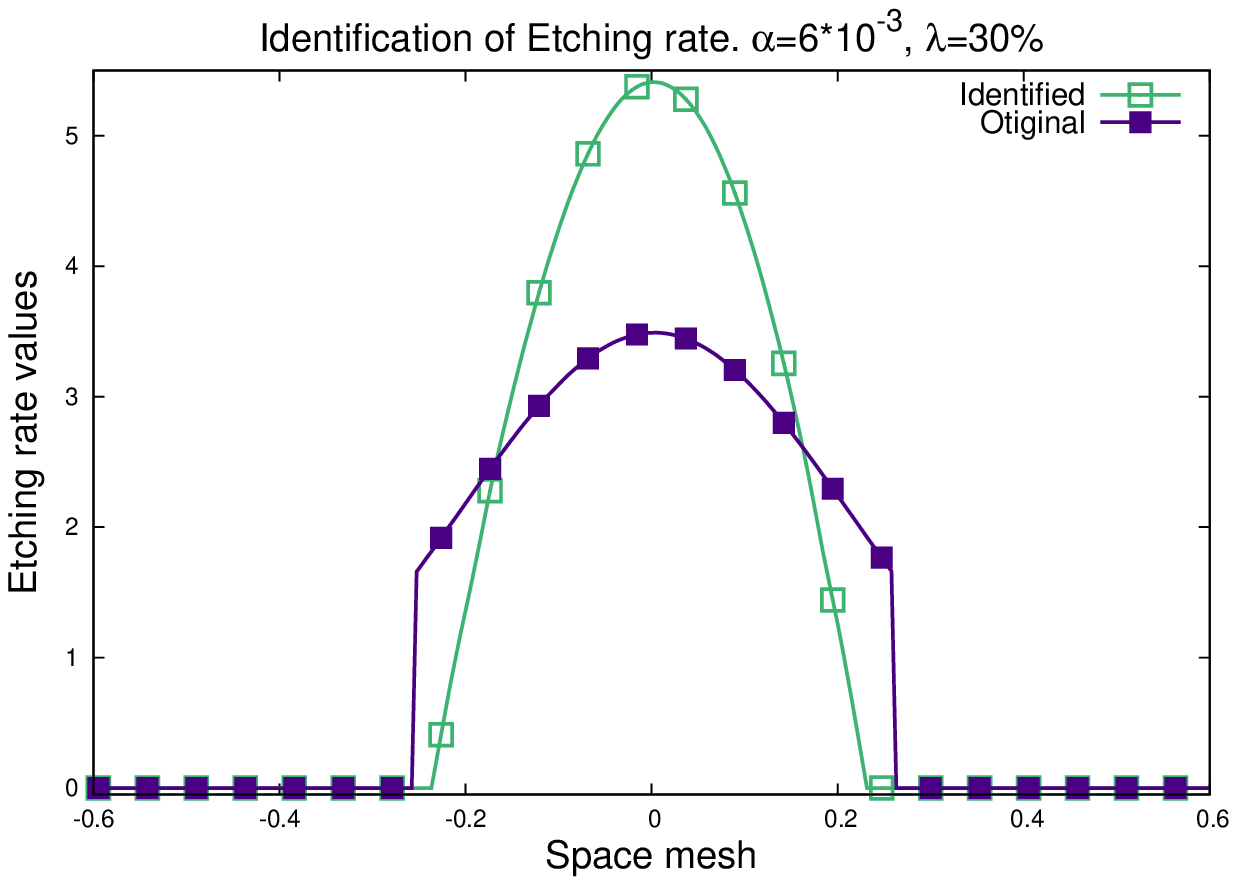}}\hspace{5pt}
				\resizebox*{6cm}{!}{\includegraphics{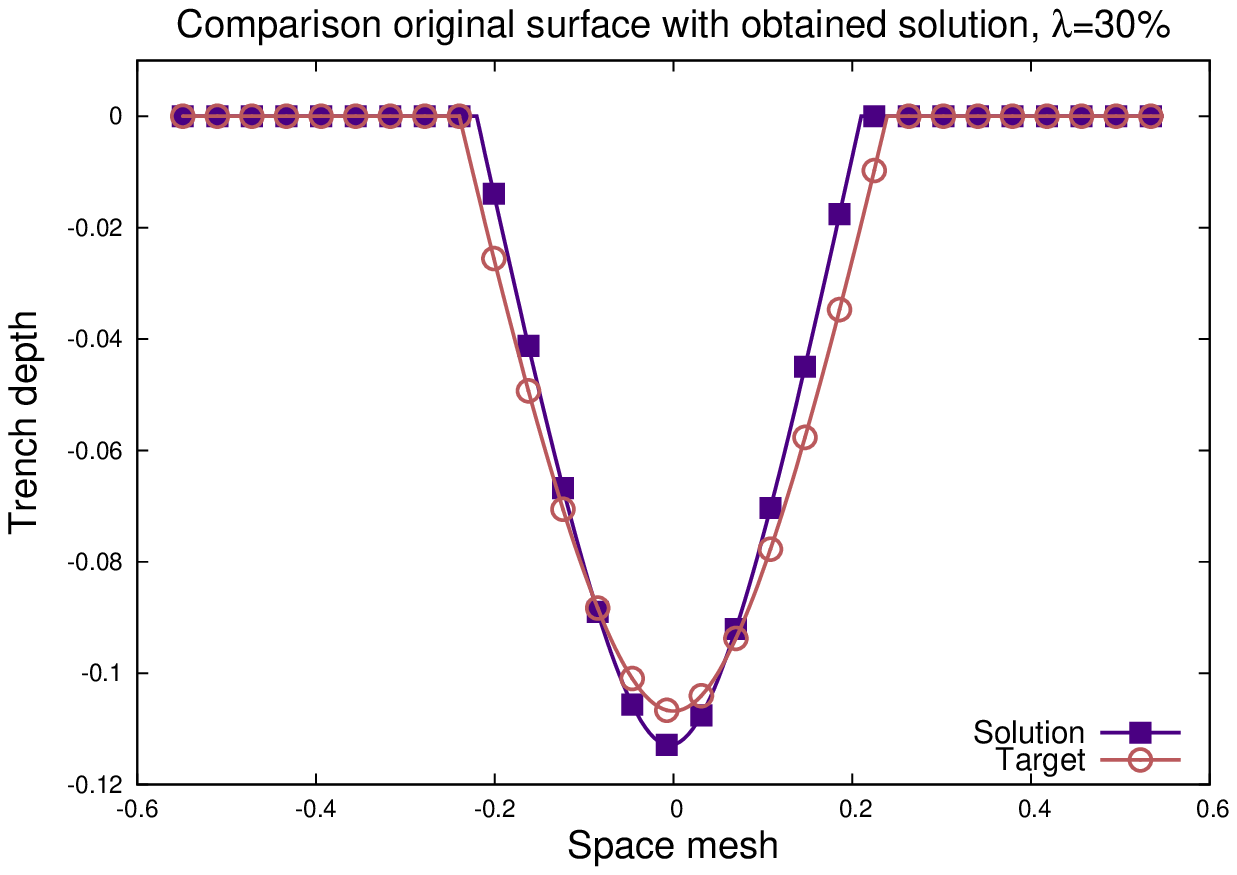}}} \\
    \caption{Results of numerical experiments in the identification the Etching rate functions for AWJM model and reconstruction of the surface shapes based on single trench profiles for measurement noises of 5\%, 15\% and 30\% respectively.
     }%
   \label{1trench_noise}
   \end{center}
\end{figure}

\begin{figure}
     \begin{center}
     \resizebox*{6cm}{!}{\includegraphics{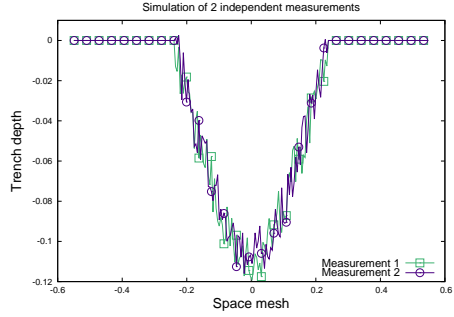}} \\
    \caption{Simulation of the 2 independent trench measurements with the level of noise = 15\%.
     }%
   \label{2trench_profile}
   \end{center}
\end{figure}

\begin{figure}
     \begin{center}
     \resizebox*{6cm}{!}{\includegraphics{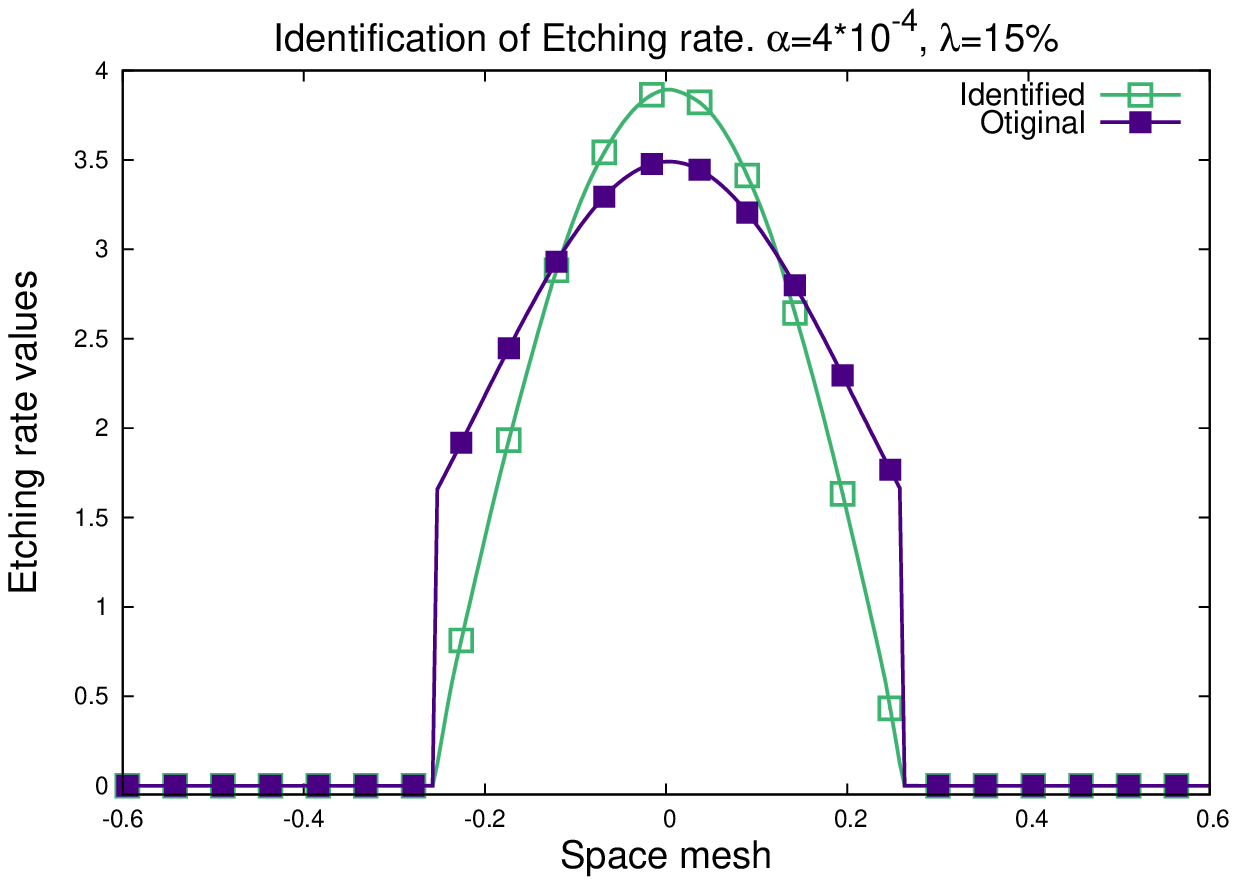}}\hspace{5pt}
     \resizebox*{6cm}{!}{\includegraphics{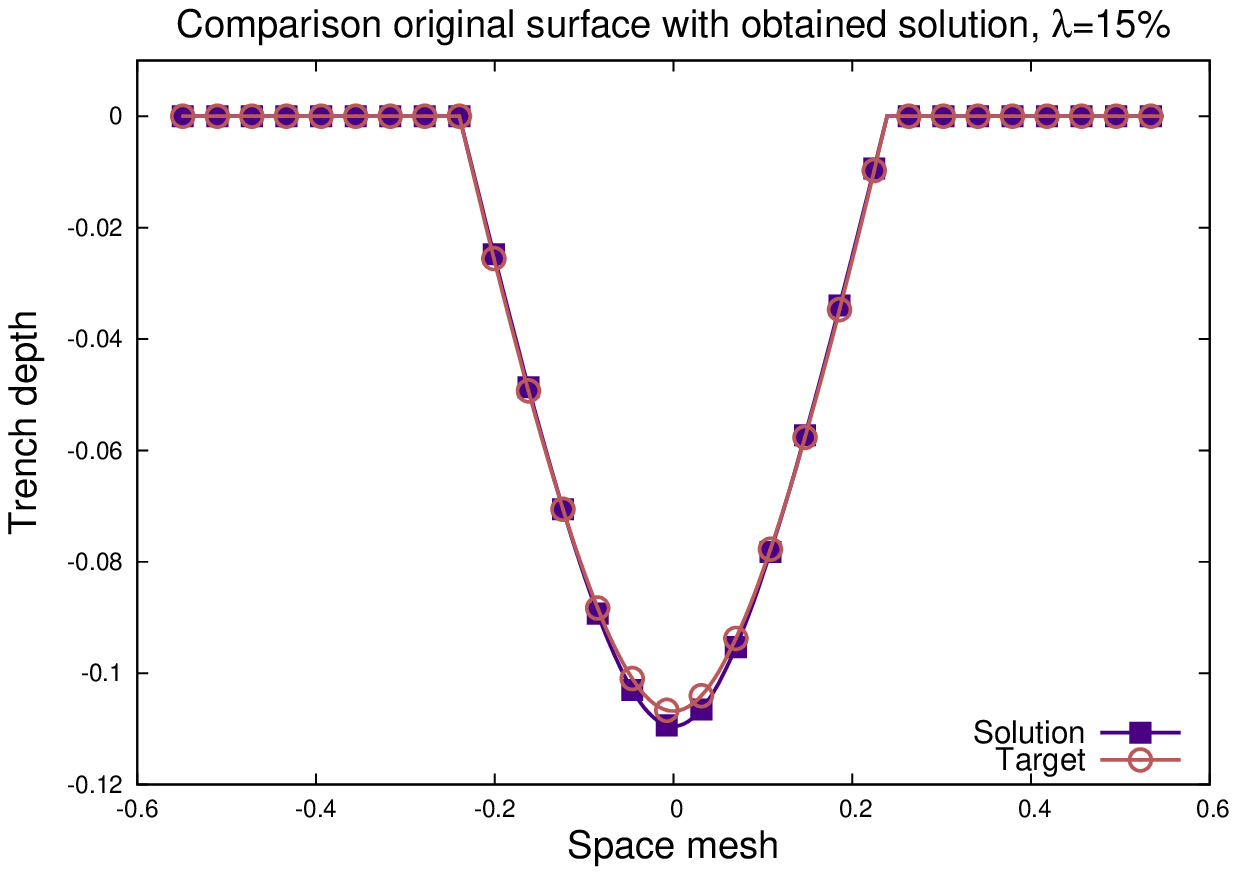}}
    \caption{Results of numerical experiments in the identification the Etching rate functions for AWJM model and reconstruction of the surface shapes based on two different measurements of one trench profile in case the level of the measurement errors is 15\%.
     }%
   \label{2tr_noise}
   \end{center}
\end{figure}

\begin{figure}
     \begin{center}
        \subfigure[Superposition of two trenches]{%
            \label{s2tr_noise}
				\resizebox*{6cm}{!}{\includegraphics{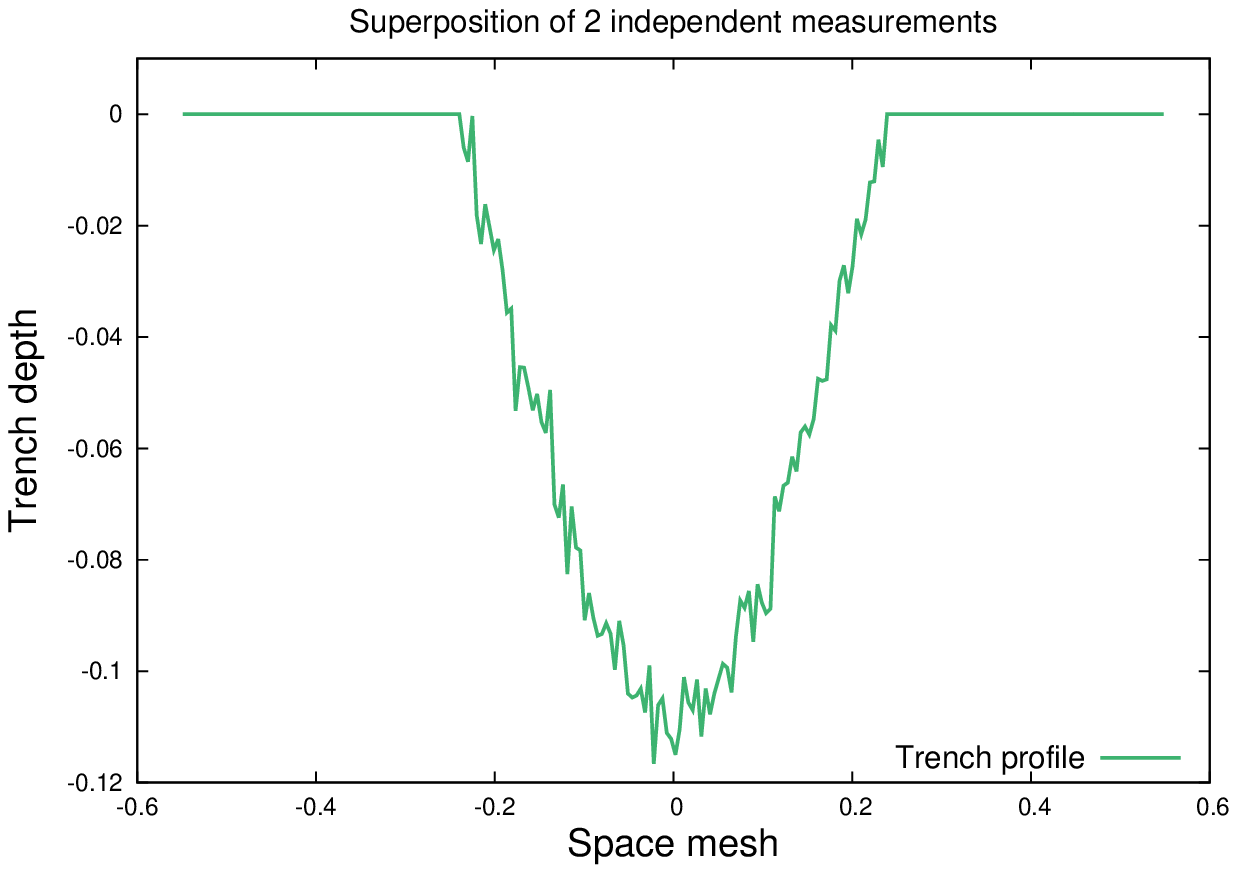}}\hspace{5pt}
				\resizebox*{6cm}{!}{\includegraphics{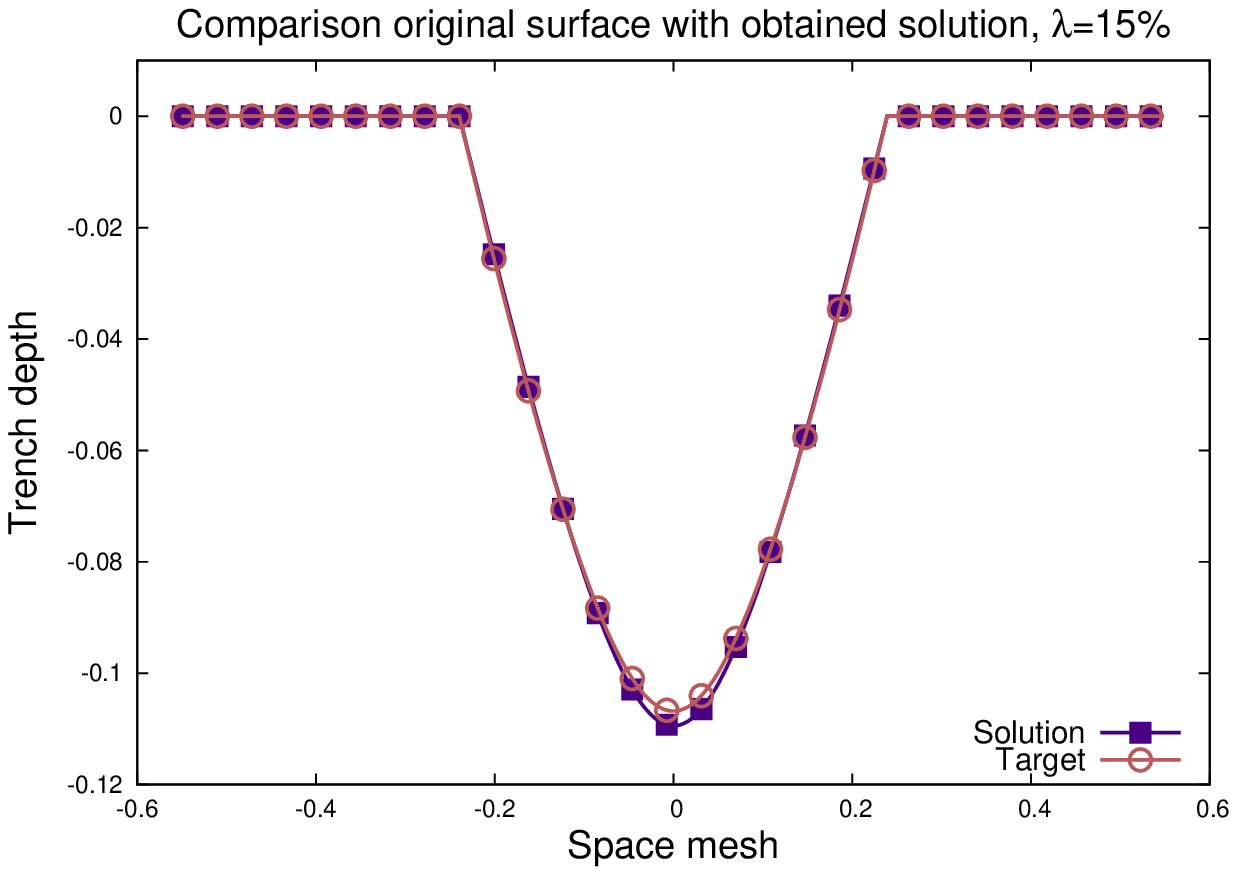}}} \\
        \subfigure[Three independent trenches]{%
            \label{3tr_noise}
				\resizebox*{6cm}{!}{\includegraphics{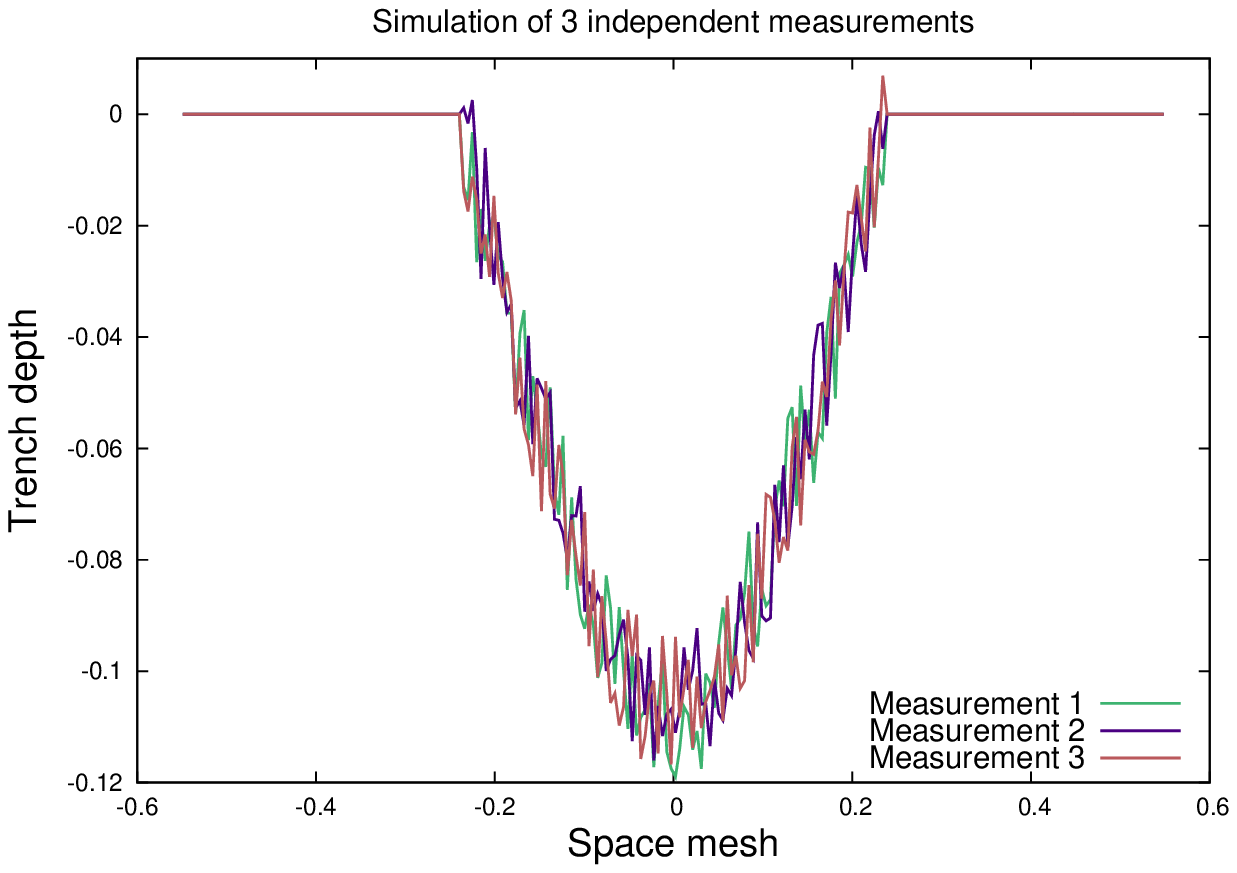}}\hspace{5pt}
				\resizebox*{6cm}{!}{\includegraphics{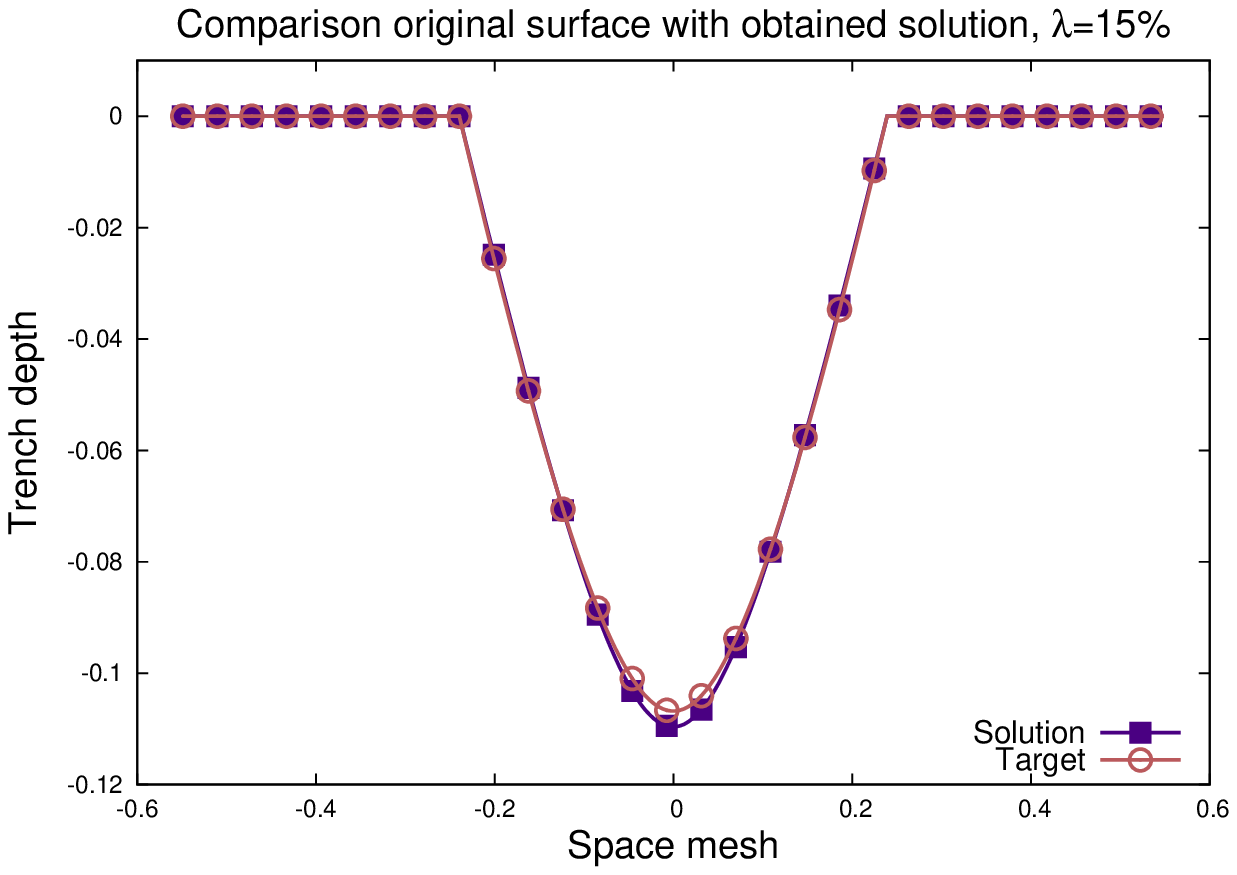}}} \\
		\subfigure[Superposition of three trenches]{%
            \label{s3tr_noise}
				\resizebox*{6cm}{!}{\includegraphics{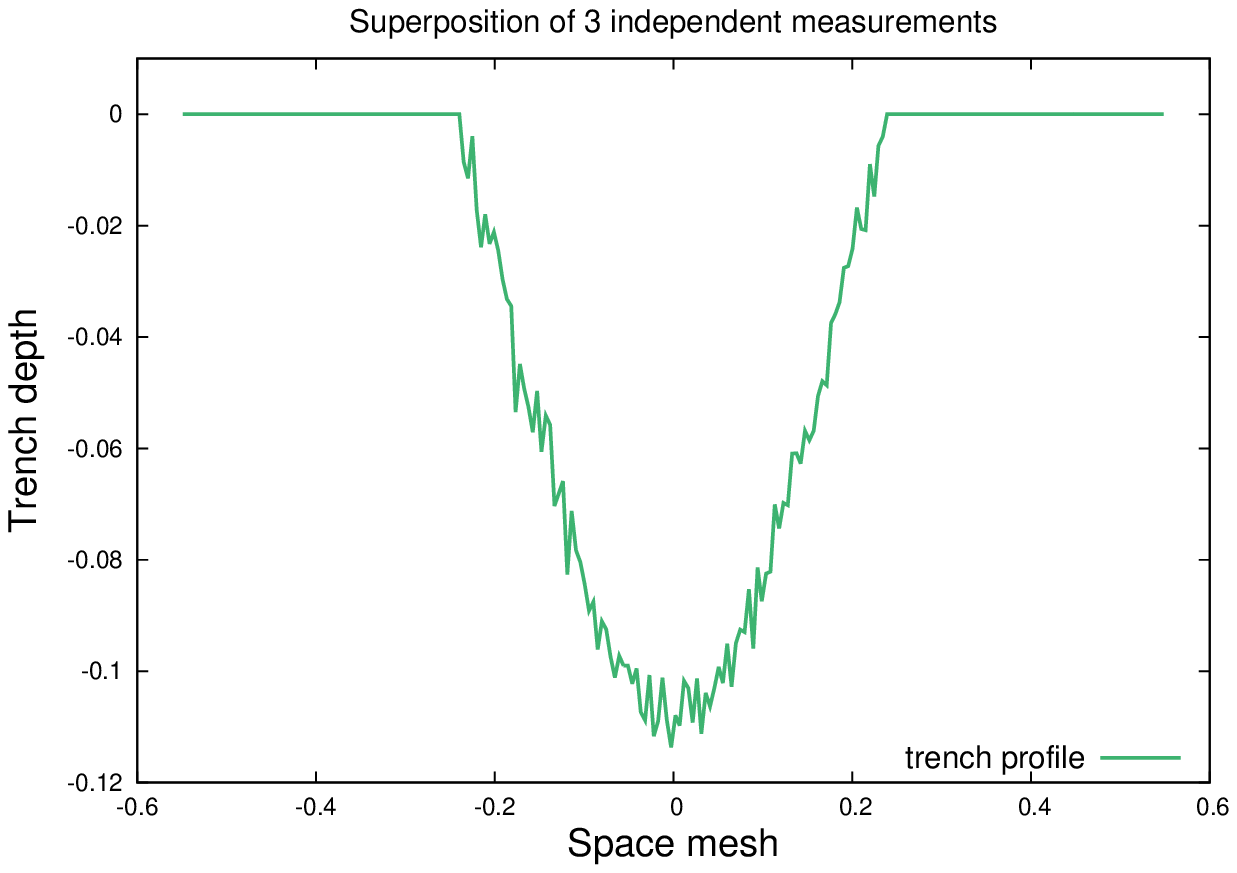}}\hspace{5pt}
				\resizebox*{6cm}{!}{\includegraphics{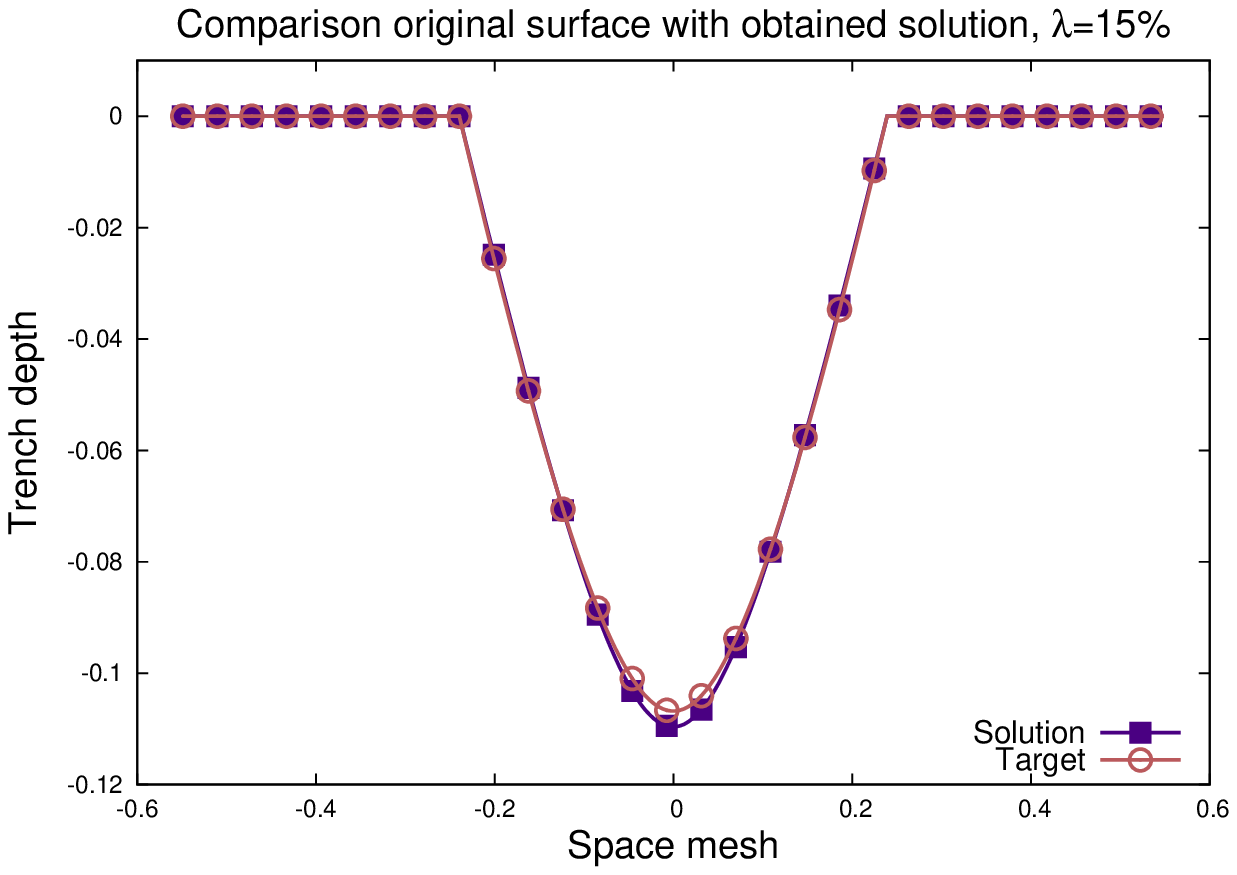}}} \\
    \caption{Results of numerical experiments in the identification the Etching rate functions for AWJM model and reconstruction of the surface shapes based on different variants of cost functions in case the level of the measurement errors is 15\%.
     }%
   \label{trenches:subfigures}
   \end{center}
\end{figure}

The identification of the model parameter $\bm E$ is based only on the experimental measurements which include the measurement errors. We simulate such data (Figure \ref{2trench_profile}) by adding a Gaussian white noise to the initial surface profile, generated by use of the Etching rate function named "Original" on Figures \ref{1trench_noise} and \ref{2tr_noise}. The goal of this study is to identify the unknown $\bm E$ the use of which in the AWJM model \eqref{2Dmodel} will form the closest trench to the initial one, named "Target" on Figures \ref{1trench_noise}, \ref{2tr_noise}, \ref{trenches:subfigures}.

In order to obtain a smooth solution, we again use the cost function \eqref{gradient cost} with the regularization term on the gradient. This will ensure the abscence of high oscillations (note that the regularization coefficient may be newly estimated through a L-curve method).

Figures \ref{1tr_5} -- \ref{1tr_30} demonstrate the results of identification of the Etching rate function $\bm E$ in case of noise with levels of 5\%, 15\% and 30\%, and comparison of reproduced trenches with original profiles.
One could also notice that even with considering a very high level of noise in the measurements, it is still possible to identify the model parameter $\bm E$. Taking these founded values, we can then model and forecast the shape of the surface, even if the form and view of the founded function $\bm E$ is not perfectly matching.

Assume now that we have two or more different measurements of exactly the same experiment. The difference between them is only in random noises that are present in measurements and influence the available input (Figure \ref{2trench_profile}).

In this case our cost function transforms to
\begin{equation}\label{gradient cost 2}
\bm J(\bm u)=\frac{1}{2}\int\limits_\Omega \left\|\bm Z(x,T)-\bm Z_{\rm exp_1}(x)\right\|^2 {\rm d}x + \frac{1}{2}\int\limits_\Omega \left\|\bm Z(x,T)-\bm Z_{\rm exp_2}(x)\right\|^2 {\rm d}x + \alpha\|\nabla \bm E\|^2,
\end{equation}
and it leads to the following results for the identification process (Figure \ref{2tr_noise}).

As it is possible to see on Figure \ref{2tr_noise} the obtained results are much more precise than in the previous case (Figure \ref{1tr_15}, only one trench profile), as the form of the identified function E becomes more symmetric and smooth. Another improvement is the symmetry of the solution obtained from very unclear measurements and increased accuracy of the trench reconstruction.

The use of more than one measurement of the only experimental data in theory will provide some kind of averaging of the trench profiles that should lead to smoothing and give higher opportunity to reconstruct the surface more precisely. Based on that assumption we introduce the superposition of the same two different measurements in the cost function as in the previous numerical case, that transforms the cost function as follows:

\begin{equation}\label{gradient cost s2}
\bm J(\bm u)=\int\limits_\Omega \left\|\bm Z(x,T)-\left( \frac{\bm Z_{\rm exp_1}(x)+\bm Z_{\rm exp_2}(x)}{2} \right) \right\|^2 {\rm d}x + \alpha\|\nabla \bm E\|^2.
\end{equation}

Numerical results of this case (Figure \ref{s2tr_noise}) show that the accuracy of the reconstruction of the trench profile is also higher in comparison with the case where only two trenches are used independently (Figure \ref{2tr_noise}).

Considering some improvment above, we now include one more additional experimental measurement, which is involved into the superposition of trenches to form the input for minimization process. First we introduce it independently from other measurements, so that the corresponding cost function is:

\begin{equation}\label{gradient cost 3}
\bm J(\bm u)=\sum_{i=1}^{3} \frac{1}{3}\int\limits_\Omega \left\|\bm Z(x,T)-\bm Z_{{\rm exp}_i}(x)\right\|^2 {\rm d}x + \alpha\|\nabla \bm E\|^2,
\end{equation}

Numerical results are displayed on the Figure \ref{3tr_noise}, where one could see that mismatch between solution and "Target" is still very low. Note that identified Etching rate function already keeps the smooth symmetric form.

Finally all the three measurements used before can also be combined together into the formation of the input for minimization process. It means that next results, presented on Figure \ref{s3tr_noise} rely on the superposition of the three different measurements of the same trench, leading to the following cost function:
\begin{equation}\label{gradient cost s3}
\bm J(\bm u)=\int\limits_\Omega \left\|\bm Z(x,T)-\left( \frac{\bm Z_{\rm exp_1}(x)+\bm Z_{\rm exp_2}(x)+\bm Z_{\rm exp_3}(x)}{3} \right) \right\|^2 {\rm d}x + \alpha\|\nabla \bm E\|^2.
\end{equation}

The difference between using two and three measurements is not very impressive due to random nature of the noise applied to the input, and could be strongly increased by involving hundreds of experimental observations to reduce the influence of the errors. But it is already possible to claim that use of three trench measurements instead of two gives more close and correct shape of the trench, and allows us to identify a more acceptable Etching rate function in the manufacturing.

\begin{figure}
     \begin{center}
     \resizebox*{9cm}{!}{\includegraphics{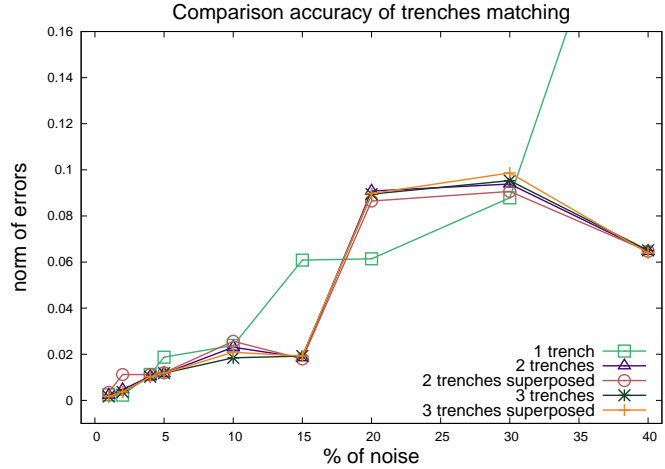}} \\
    \caption{Accuracy of matching the obtained solutions to original profiles for different cost functions considering different levels of measurement errors.
     }%
   \label{propagation}
   \end{center}
\end{figure}

Figure \ref{propagation} gives a good overview of the improvement in the identification achieved by introducing different approaches to define the cost function. Keeping the acceptable level of errors (less than 10\%) in matching between experimental measurement and obtained solution, introduction of more complex cost functions \eqref{gradient cost 2} and \eqref{gradient cost 3} gives us the opportunity to identify more applicable and convenient AWJM model parameters. In general the more trenches measurements are available the higher accuracy in the surface prediction can be achieved due to elimination of the measurement errors by averaging the input data.

Cost functions \eqref{gradient cost 2} and \eqref{gradient cost s2} are not identical, but they theoretically have the same gradient. However numerical implementation gives different results and flexibility to find more suitable realization for each particular problem, what has to be taken into account. The same situation is for the pair of cost functions \eqref{gradient cost 3} and \eqref{gradient cost s3}. More detailed results of comparison between all proposed approaches are presented in Table \ref{tablichka}.

\begin{table}
\tbl{Comparison of the accuracy in the prediction of the trench profile, corresponding to different cases of the cost functions and different levels of applied noise.}
{\begin{tabular}[l]{@{}lcccccccc}\toprule
  Trenches & 1\% & 2\% & 5\% & 10\% & 15\% & 20\% & 30\% & 40\% \\
\colrule
  Single & $2.37\times 10^{-3}$ & $2.24\times 10^{-3}$ & $1.87\times 10^{-2}$ & $2.38\times 10^{-2}$ & $6.08\times 10^{-2}$ & $6.14\times 10^{-2}$ & $8.78\times 10^{-2}$ & 0.272735 \\
\colrule
  2 independent & $2.23\times 10^{-3}$ & $4.66\times 10^{-3}$ & $1.15\times 10^{-2}$ & $2.30\times 10^{-2}$ & $1.84\times 10^{-2}$ & $9.08\times 10^{-2}$ & $9.38\times 10^{-2}$ & $6.46\times 10^{-2}$ \\
\colrule  
  2 superposed  & $3.43\times 10^{-3}$ & $1.12\times 10^{-2}$ & $1.19\times 10^{-2}$ & $2.56\times 10^{-2}$ & $1.79\times 10^{-2}$ & $8.64\times 10^{-2}$ & $9.06\times 10^{-2}$ & $6.43\times 10^{-2}$ \\
\colrule  
  3 independent & $1.68\times 10^{-3}$ & $3.66\times 10^{-3}$ & $1.18\times 10^{-2}$ & $1.84\times 10^{-2}$ & $1.92\times 10^{-2}$ & $8.94\times 10^{-2}$ & $9.53\times 10^{-2}$ & $6.50\times 10^{-2}$ \\
\colrule  
  3 superposed  & $1.68\times 10^{-3}$ & $3.65\times 10^{-3}$ & $1.18\times 10^{-2}$ & $2.08\times 10^{-2}$ & $1.92\times 10^{-2}$ & $8.97\times 10^{-2}$ & $9.86\times 10^{-2}$ & $6.39\times 10^{-2}$ \\
\botrule
\end{tabular}}

\label{tablichka}
\end{table}

One could notice that in most of the results the use of several trenches instead of only one can improve the accuracy in the parameters identification, leading to reducing the errors in the prediction of the surface profile up to 20\% in cases of low level of noise. In cases of very noisy input data (see columns related to noise higher than 20\%) the use of several measurements plays strong role in decreasing the mismatch in reconstruction of the trenches. This effect could be explained by nature of the distribution of the applied noise. Certainly, it should be noted that sometimes only one measurement is available, and it might be enough to obtain the model parameters required to reconstruct the profile. Moreover, superposition of the trenches involved in the identification process which can be interpreted as an average of the input experiments usually also leads to higher accuracy, but mostly on the "long distance" -- the quantity of available measurements.

\section{Conclusion}
\label{sec_concl}

The identification of model parameters especially from the noisy data is a challenging problem because of its ill-posedness. This paper has presented the application of inverse problems theory, based on minimization problems in the real production. The identification of optimal model parameters indeed gives a chance to model and predict the trench profile for AWJ machining. We illustrated how an even small level of errors or noise can influence the identification results and lead to significant errors in the determination. Also we explained how the high level of noise could completely change the identification process and why it is necessary to keep the regularization terms to get a more precise structure and shape of the reconstructed trench.

We presented a way to estimate the required surface profile in the lack of knowledge of exact AWJ model parameters. Proposed AWJM model coupled with measurement errors which were taken into account can precisely model the surface shape by identifying optimal values of the model parameters even with very poor and inaccurate experimental observations.

In this paper we considered measurement errors that were included in the mathematical model, which we tried to compensate by adding Tikhonov regularization terms. Development of this study in the case of nonsteady AWJM process in 3D, including the variation of the jet feed speed, leads to complication of the model and numerical realization of the identification process. All these aspects may be studied more precisely in a future work.

\section*{Acknowledgements}

The authors would like to acknowledge the funding support of the EU-FP7-ITN (Grant no. 316560) for the works presented as a part of the STEEP ITN project. The authors thank Mr. Pablo Lozano Torrubia from the University of Nottingham for his help in the experimental work.

\bibliographystyle{gIPE}
\bibliography{GROZA}

\begin{thebibliography}{10}
\providecommand{\url}[1]{\normalfont{#1}}
\providecommand{\urlprefix}{Available from: }

\bibitem{Lavr1}
Lavrentiev~MM, Romanov~VG, Shishatskii~SP. Ill-posed problems of mathematical
  problems. Vol.~64. Providence: AMS; 1986.

\bibitem{Tarantola}
Tarantola~A. Inverse problem theory and methods for model parameter estimation.
  Philadelphia: SIAM; 2005.

\bibitem{Aster}
Aster~RC, Borchers~B, Thurber~CH. Parameter estimation and inverse problems
  (second edition). Elsevier; 2013.

\bibitem{Tikh_Arsenin}
Tikhonov~AN, Arsenin~VY. Solutions of ill-posed problems. New York: John Wiley
  \& Sons; 1977.

\bibitem{Tikh_Glas}
Tikhonov~AN, Glasko~VB. Use of the regularization method in non-linear
  problems. Zh Vychisl Mat Mat Fiz. 1965;\hspace{0pt}5(3):463--473.

\bibitem{Tautenhahn2}
Tautenhahn~U. Lavrentiev regularization of nonlinear ill-posed problems.
  Vietnam Journal of Mathematics. 2004;\hspace{0pt}32:29--41.

\bibitem{Denisov}
Denisov~AM, Zaharov~EV, Kalinin~AV, Kalinin~VV. Application of tikhonov
  regularization method for numerical solution of inverse problem of
  electrocardiography. MSU Vestnik. 2008;\hspace{0pt}15.

\bibitem{Barbara2}
Kaltenbacher~B, Neubauer~A, Scherzer~O. Iterative regularization methods for
  nonlinear ill-posed problems. Vol.~6 of Radon Series on Computational and
  Applied Mathematics. Berlin: Walter de Gruyter GmbH \& Co; 2008.

\bibitem{Heinz_Hanke}
Engl~HW, Hanke~M, Neubauer~A. Regularization of inverse problems. Vol. 375 of
  Mathematics and its Applications. Dordrecht: Kluwer Academic Publishers
  Group; 1996.

\bibitem{Dragos1}
Axinte~DA, Srinivasu~DS, Billingham~J, Cooper~M. Geometrical modelling of
  abrasive waterjet footprints: A study for $90^{\circ}$ jet impact angle. CIRP
  Annals -- Manufacturing Technology. 2010;\hspace{0pt}59:341--346.

\bibitem{Dragos2}
Kong~MC, Anwar~S, Billingham~J, Axinte~DA. Mathematical modelling of abrasive
  waterjet footprints for arbitrarily moving jets: Parti -- single straight
  paths. International Journal of Machine Tools \& Manufacture.
  2012;\hspace{0pt}53:58--68.

\bibitem{Dragos3}
Billingham~J, Miron~CB, Axinte~DA, Kong~MC. Mathematical modelling of abrasive
  waterjet footprints for arbitrarily moving jets: Partii -- overlapped single
  and multiple straight paths. International Journal of Machine Tools \&
  Manufacture. 2013;\hspace{0pt}68:30--39.

\bibitem{Pablo}
Lozano~Torrubia~P, Billingham~J, Axinte~DA. Stochastic simplified modelling of
  abrasive waterjet footprints. Proc R Soc A. 2016;\hspace{0pt}:2016 472
  20150836.

\bibitem{Bilbao}
Bilbao~Guillerna~A, Axinte~DA, Billingham~J. The linear inverse problem in
  energy beam processing with an application to abrasive waterjet machining.
  International Journal of Machine Tools \& Manufacture.
  2015;\hspace{0pt}99:34--42.

\bibitem{Gil89}
Gilbert~JC, Lemar{\'e}chal~C. Some numerical experiments with variable storage
  quasi-newton algorithms. Math Prog. 1989;\hspace{0pt}45:407--435.

\bibitem{Liu89}
Liu~DC, Nocedal~J. On the limited memory bfgs method for large scale
  optimization. Math Prog. 1989;\hspace{0pt}45:503--528.

\bibitem{Didier1}
Veers{\'e}~F, Auroux~D. Some numerical experiments on scaling and updating
  l-bfgs dianogal preconditioners. INRIA; 2000. Research Report 3858.

\bibitem{N2QN1}
Gilbert~JC, Lemar{\'e}chal~C. 2013;
  \urlprefix\url{https://who.rocq.inria.fr/Jean-Charles.Gilbert/modulopt}.

\bibitem{TAPENADE}
Hasco{\"e}t~L, Pascual~V. The tapenade automatic differentiation tool:
  Principles, model, and specification. ACM Transactions On Mathematical
  Software. 2013;\hspace{0pt}39(3).

\bibitem{Quarteroni}
Quarteroni~A. Numerical models for differential problems. Springer--Verlag;
  2009.

\bibitem{Lions}
Lions~JL. Optimal control of systems governed by partial differential
  equations. Springer--Verlag; 1971.

\bibitem{Laws_Hanson}
Lawson~CL, Hanson~RJ. Solving least squares problems. Philadelphia (PA): SIAM;
  1995.

\bibitem{Hensen_Oleary}
Hansen~PC, O'Leary~DP. The use of the l-curve in the regularization of discrete
  ill-posed problems. SIAM J Sci Comput. 1993;\hspace{0pt}14:1487--1503.

\end{thebibliography}

\end{document}